\documentclass[11pt,a4paper]{article}
\usepackage{amsfonts,amsmath,amssymb,graphics,amsthm,amscd}
\usepackage[french]{babel}
\usepackage[T1]{fontenc}
\usepackage[latin1]{inputenc}

\newcommand{\N}{\mathbb{N}}
\newcommand{\Z}{\mathbb{Z}}

\def\F{\mbox{\it I\hspace{-.15em}F}}

\def\G{\Gamma }
\def\a{\alpha }

\def\s{\sigma }

\def\p{\varphi }

\def\t{\theta }

\newcommand{\fonc}[5]{#1 : \left\{ \begin{array}{c c c}
                            #2 & \longrightarrow & #3 \\ 
                            #4 & \longmapsto & #5 
                            \end{array} \right.}

\newcommand{\bij}[5]{#1 : \left\{ \begin{array}{c c c}
                            #2 & \stackrel{\sim}{\longrightarrow} & #3 \\ 
                            #4 & \longmapsto & #5 
                            \end{array} \right.}

\title{SUR LES AUTOMORPHISMES ET LA RIGIDIT\'E DES GROUPES DE COXETER \`A ANGLES DROITS.}
\author{Anatole CASTELLA}
\date{12 Décembre 2004}

\begin{document}

\maketitle

\theoremstyle{definition}
\newtheorem{defi}{Définition}
\newtheorem{defis}{Définitions}
\newtheorem*{defis*}{Définitions}
\newtheorem*{defi*}{Définition}
\newtheorem{rema}{Remarques}
\newtheorem{nota}{Notation}
\newtheorem{ex}{Exemple}
\newtheorem{propr}{Propriété}
\newtheorem{comms}{Commentaires}
\theoremstyle{plain}

\newtheorem{propo}{Proposition}
\newtheorem*{propo*}{Proposition}
\newtheorem{lem}{Lemme}
\newtheorem{theo}{Théorème}
\newtheorem*{theo*}{Théorème}
\newtheorem{cor}{Corollaire}
\newtheorem*{cor*}{Corollaire}

\section*{Introduction.}

Un {\it groupe de Coxeter} est un groupe $W$ admettant un sous-ensemble générateur $S$ formé d'éléments d'ordre 2 tel que, si l'on note, pour $s, t \in S$, $m_{s,t}$ l'ordre du produit $st$ dans $W$, le groupe $W$ soit donné par la présentation : 
$$W = <S \mid (st)^{m_{s,t}} = 1, \textrm{ si } s,t \in S  \textrm{ avec }  m_{s,t}\not= \infty>.$$

On dit alors que $(W,S)$ est un {\it système de Coxeter} (admis par $W$), que l'ensemble $S$ est un {\it ensemble de Coxeter} pour $W$ et que le cardinal de $S$ est le {\it rang} du système $(W,S)$. On note ${\cal S}(W)$ l'ensemble des ensembles de Coxeter pour $W$.

Généralement, on se donne un système de Coxeter $(W,S)$ par la matrice $\G_{(W,S)} = (m_{s,t})_{s,t\in S}$, que l'on représente par un graphe d'ensemble de sommets $S$, avec une arête étiquetée $m_{s,t}$ entre les sommets $s$ et $t$ lorsque $m_{s,t} \geq 3$. On dit que la matrice $\G_{(W,S)}$ est {\it le type} de $(W,S)$ (et est {\it un type} admis par $W$). C'est une {\it matrice de Coxeter}, c'est-à-dire une matrice symétrique à coefficients dans $\N^\star \cup \{\infty\}$, avec des 1 sur la diagonale et uniquement sur la diagonale. On peut montrer (cf. [B], Ch. V, \S \ 4, n° 3) que toute matrice de Coxeter est, à isomorphisme de matrices près, le type d'un système de Coxeter (où l'on appelle {\it isomorphisme de matrices} de $(m_{s,t})_{s,t\in S}$ sur $(m'_{s',t'})_{s',t'\in S'}$ toute bijection de $S$ sur $S'$ qui respecte les coefficients). \\

Soit $W$ un groupe de Coxeter. 

Si $S\in {\cal S}(W)$ et si $\a$ est un automorphisme de $W$, alors $\a(S)\in {\cal S}(W)$ et les types $\G_{(W,S)}$ et $\G_{(W,\a(S))}$ sont isomorphes (via $S \rightarrow \a(S)$, $s \mapsto \a(s)$). Inversement, pour $S,~S' \in {\cal S}(W)$, tout isomorphisme de $\G_{(W,S)}$ sur $\G_{(W,S')}$ se prolonge en un automorphisme de $W$ qui envoie $S$ sur $S'$. On est alors amené à considérer les deux problèmes suivants : 
\begin{enumerate}
\item Déterminer le groupe $Aut(W)$.
\item \label{classes d'iso}Déterminer les classes d'isomorphisme des types admis par $W$ ou, de façon équivalente, les orbites de l'action de $Aut(W)$ sur ${\cal S}(W)$ : $Aut(W)\times {\cal S}(W)\rightarrow {\cal S}(W)$, $(\a,S)\mapsto \a(S)$.
\end{enumerate}
On dit qu'un groupe de Coxeter est {\it rigide} lorsqu'il admet un unique type, à isomorphisme près, c'est-à-dire lorsque l'action décrite en \ref{classes d'iso} ci-dessus est transitive (voir aussi la définition \ref{defs rigidité} en section \ref{Conclusion rigidité.} ci-dessous).

De nombreux résultats existent sur ces deux problèmes, pour différentes classes de groupes de Coxeter, et en particulier pour la classe des groupes de Coxeter "à angles droits" : \\

\begin{defi}[Angles droits]\label{def angles droits}Soit $(W,S)$ un système de Coxeter de type $\G$. On dit que $(W,S)$ et $\G$ sont {\it à angles droits} si $\G$ est à coefficients dans $\{1,2,\infty\}$. 

On montre que, si un groupe de Coxeter $W$ admet un système de Coxeter à angles droits, alors tous les systèmes de Coxeter qu'il admet sont à angles droits et ont même rang (voir la section \ref{Le cas à angles droits.}). On dit alors que le {\it groupe} de Coxeter $W$ est à angles droits et que le rang commun des systèmes de Coxeter qu'il admet est le rang du {\it groupe} $W$.\\
\end{defi} 

Soit $W$ un groupe de Coxeter à angles droits (de rang fini ou infini). 

Notons $\pi : W \rightarrow W^{ab}, \ w \mapsto \overline w$ le morphisme canonique de $W$ sur son abélianisé. Comme le sous-groupe dérivé est un sous-groupe caractéristique, $\pi$ induit le morphisme $\pi_{Aut} : Aut(W) \longrightarrow Aut(W^{ab}), \ \a \longmapsto (\overline \a : \overline w \mapsto \overline {\a(w)}\, )$. Notons $F$ l'ensemble des éléments d'ordre fini de $W$, et posons $\overline F = \pi(F)$. Comme les automorphismes de $W$ stabilisent $F$, $\pi_{Aut}$ est à valeurs dans le sous-groupe $Aut(W^{ab},\overline F)$ de $Aut(W^{ab})$ constitué des automorphismes de $W^{ab}$ qui stabilisent $\overline F$. En 1988, J. Tits a établi dans [T] (Corollaire 1) que la suite :
$$\{1\}\hookrightarrow \ker(\pi_{Aut}) \hookrightarrow Aut(W) \stackrel{\pi_{Aut}}{\longrightarrow}Aut(W^{ab},\overline F) \rightarrow \{1\}$$
est exacte et scindée. Il a ensuite entamé l'étude de $\ker(\pi_{Aut})$ (qu'il note $Aut^\circ(W)$), fourni des pistes à suivre pour étudier $Aut(W^{ab},\overline F)$ et décrit ce dernier groupe dans quelques exemples. Dix ans plus tard, B. Mühlherr a donné dans [M] une présentation par générateurs et relations de $\ker(\pi_{Aut})$ (qu'il note $Spe(W)$), dans le cas où $W$ est de rang fini.

Au début des années 2000, D.G. Radcliffe et T. Hosaka ont montré respectivement dans [R] et dans [H] que les groupes de Coxeter à angles droits et de rang fini sont rigides. \\

Le but de cet article est de donner une démonstration de ce dernier résultat qui ne fait pas intervenir la finitude du rang (donc de montrer que {\it tous} les groupes de Coxeter à angles droits sont rigides)\footnote{Après la première diffusion de cet article, D.G. Radcliffe m'a signalé que cette généralisation est également conséquence d'un résultat de sa thèse de doctorat ([R2], Ch. 5, Théorème principal).} et de décrire le groupe $Aut(W^{ab},\overline F)$ pour une vaste classe de groupes de Coxeter à angles droits (contenant les groupes de Coxeter à angles droits et de rang fini), à savoir ceux que nous appelons "d'épaisseur finie" (section \ref{profondeur}, définition \ref{def profondeur finie}).\\

La première partie est consacrée à quelques rappels sur les groupes et sur les groupes de Coxeter, qui nous permettent de fixer les notations utilisées dans les parties suivantes.

Dans la seconde partie, nous présentons des résultats généraux sur les ensembles munis d'une relation binaire symétrique et réflexive (dont la motivation est expliquée ci-dessous). La terminologie employée est la suivante : \\

\begin{defi*}[Section \ref{partie outils}, définitions \ref{def relations}, \ref{defs parties commutatives}, \ref{defs bicommutants cellules} et \ref{def noyaux}]Nous appelons {\it ensemble à relation} tout couple $(E,R)$ où $E$ est un ensemble et $R$ une relation binaire symétrique et réflexive sur $E$. 

Soit $(E,R)$ un ensemble à relation. Nous notons ${\cal P}(E)$ l'ensemble des parties de $E$ et $[E]$ l'ensemble des parties finies de $E$. Ce sont des $\F_2$-espaces vectoriels (cf. section \ref{Généralités sur les groupes.}). Nous notons $[E]_c$ l'ensemble des parties finies de $E$ constituées d'éléments deux à deux en relation $R$.

Soit $C : {\cal P}(E) \rightarrow {\cal P}(E)$, $X \longmapsto \{y\in E \mid \forall x\in X, \ yRx \}$. Nous appelons {\it cellules} les images par $C^2$ des singletons de $E$ et {\it noyaux} les classes d'équivalence de la relation (d'équivalence) $C(\{x\}) = C(\{y\})$ sur $E$. Nous notons ${\cal C}(E,R)$ l'ensemble des cellules et ${\cal N}(E,R)$ l'ensemble des noyaux. Pour $x \in E$, nous disons que la cellule $C^2(\{x\})$ est la cellule de (ou définie par) $x$ et nous notons $N(x)$ le noyau (la classe d'équivalence) de $x$.

Si $(E,R)$ et $(E',R')$ sont deux ensembles à relation, nous appelons {\it isomorphisme (d'ensembles à relation)} de $(E,R)$ sur $(E',R')$ toute bijection $f$ de $E$ sur $E'$ telle que $(f\times f)(R) = R'$.\\
\end{defi*}

Si $\G = \G_{(W,S)}$ est un type admis par $W$, nous munisons $S$ de la relation $R_\G$ de commutation dans $S$. Nous obtenons ainsi un ensemble à relation $(S,R_\G)$ qui détermine entièrement $\G$ (qui est à angles droits), puisque l'on a, pour $s,~t \in S$ distincts, $m_{s,t} = 2 \Leftrightarrow sR_\G t$ et $ m_{s,t} = \infty \Leftrightarrow (s,t)\not \in R_\G$ (notons que $[S]_c$ est alors l'ensemble des parties commutatives et finies de $S$). 

Pour $S,~S'\in {\cal S}(W)$ et $\G = \G_{(W,S)}$, $\G' = \G_{(W,S')}$, on voit que les notions d'isomorphisme de matrices de $\G$ sur $\G'$ et d'isomorphisme d'ensembles à relation de $(S,R_{\G})$ sur $(S',R_{\G'})$ coïncident. 

De plus, à la suite de J. Tits, nous remarquons qu'il existe un isomorphisme ($\F_2$-linéaire) de $[S]$ sur $[S']$ envoyant $[S]_c$ sur $[S']_c$ (voir la démonstration du corollaire 1 de [T], où $[S]_c$ est noté $F(\G)$, ou la remarque \ref{relations et singletons} a ci-dessous). 

Pour démontrer que les groupes de Coxeter à angles droits (de rang quelconque) sont rigides, on voit donc qu'il suffit de montrer le résultat suivant : \\

\begin{theo*}[Section \ref{sur les isomorphismes ...}, assertion 2 du théorème \ref{theorème principal}]Soient $(E,R)$ et $(E',R')$ deux ensembles à relation. S'il existe un isomorphisme ($\F_2$-linéaire) de $[E]$ sur $[E']$ envoyant $[E]_c$ sur $[E']_c$, alors les ensembles à relation $(E,R)$ et $(E',R')$ sont isomorphes.\\
\end{theo*}

Dans la troisième partie, nous examinons les notions de rigidité et de rigidité forte pour les groupes de Coxeter à angles droits  (cf. définitions \ref{defs rigidité} et \ref{def rigidité forte}). Nous obtenons en particulier, comme corollaire immédiat du Théorème \ref{theorème principal}, le résultat suivant : \\

\begin{theo*}[Section \ref{Conclusion rigidité.}, Théorème \ref{Les groupes de Coxeter à angles droits sont rigides}]Les groupes de Coxeter à angles droits sont rigides.\\
\end{theo*}

Dans la quatrième partie, nous étudions le groupe $Aut(W^{ab},\overline F)$. Comme l'avait remarqué J. Tits dans [T], pour tout $S\in {\cal S}(W)$, nous pouvons identifier $Aut(W^{ab},\overline F)$ au groupe $Aut([S],[S]_c)$ constitué des automorphismes de $[S]$ stabilisant $[S]_c$ (voir la section \ref{Le cas à angles droits.} ci-dessous). Fixant $S \in {\cal S}(W)$, c'est ce dernier groupe $Aut([S],[S]_c)$ que nous décrivons ici. L'étude que nous effectuons est inspirée des idées exposées dans la partie finale de [T].

Posons $\G = \G_{(W,S)}$, ${\cal C} = {\cal C}(S,R_\G)$ et ${\cal N} = {\cal N}(S,R_\G)$. La proposition \ref{sous-espaces [T] permutés} de la partie \ref{sur les isomorphismes ...} montre en particulier que le groupe $Aut([S],[S]_c)$ agit sur l'ensemble $\{[T] \mid T \in {\cal C}\}$, via $(\p,[T])\mapsto \p([T])$. Notant $K(\G)$ le noyau de cette action, nous exhibons un sous-groupe $G$ du groupe $Aut(\G)$ des automorphismes de la matrice $\G$ tel que l'on ait $Aut([S],[S]_c) = K(\G) \rtimes G$ (cf. section \ref{premier dévissage}). 

Nous donnons ensuite une caractérisation simple des éléments de $K(\G)$ parmi les éléments de $Aut([S])$ (section \ref{le groupe K}, proposition \ref{caract de F et K}) et montrons que le groupe $K(\G)$ se décompose lui-même en le produit semi-direct $K(\G) = K^\circ(\G) \rtimes D(\G)$ (section \ref{le groupe K}, proposition \ref{decomposition K(Gamma)}), où $D(\G)$ désigne le groupe constitué des automorphismes de $[S]$ qui, pour tout $N \in {\cal N}$, stabilisent le sous-espace $[N]$ de $[S]$ (remarquons que, ${\cal N}$ étant une partition de $S$, on a $[S] = \oplus_{N\in {\cal N}}[N]$ et $D(\G) \approx \prod_{N\in {\cal N}}Aut([N])$) et où $K^\circ(\G)$ est le sous-groupe de $K(\G)$ suivant : \\

\begin{defi*}[Section \ref{le groupe K}, notation \ref{Le groupe Krond}]Soit $s\in S$. Posons $(C^2(\{s\}))^\star = C^2(\{s\}) \setminus N(s)$. 

Alors $K^\circ(\G) = \{\p \in K(\G) \mid \forall ~s\in S, \ \p(\{s\}) \in \{s\} + [(C^2(\{s\}))^\star]\}$.\\
\end{defi*}

La section \ref{K0} est consacrée à l'étude du groupe $K^\circ(\G)$. La méthode consiste essentiellement à déterminer une partition de $S$ en parties "d'épaisseur nulle" (section \ref{Parties de profondeur nulle.}, définition \ref{def parties de profondeur nulle}), qui se traduit par une décomposition de $K^\circ(\G)$ en produits semi-directs itérés de certains de ses sous-groupes, simples à décrire (section \ref{Parties de profondeur nulle.}, formules (E1) ou (E2), et section \ref{profondeur}, théorème \ref{décomp de KY}). C'est à cette étape qu'il faut supposer $S$ "d'épaisseur finie" (section \ref{profondeur}, définition \ref{def profondeur finie}), afin que les produits semi-directs itérés mentionnés soient en nombre fini.

Les résultats obtenus sont illustrés par quelques exemples à la fin de cet article.

\section{Préliminaires.}

\subsection{Généralités sur les groupes.}\label{Généralités sur les groupes.}

Soit $W$ un groupe. Nous notons $F(W)$, ou plus simplement $F$, l'ensemble des éléments d'ordre fini de $W$. Soit $\pi : W\rightarrow W^{ab}, \ w \mapsto \overline w$ le morphisme canonique de $W$ sur son abélianisé. Nous notons $\overline F = \overline {F(W)}$ l'image de $F(W)$ par $\pi$.

Nous notons $Aut(W)$ le groupe des automorphismes de $W$ et $Int(W)$ le groupe de ses automorphismes intérieurs. Pour $X\subseteq W$, notons $Aut(W,X) = \{\a \in Aut(W) \mid \a(X) = X\}$ le sous-groupe de $Aut(W)$ formé des automorphismes de $W$ qui stabilisent $X$. \\

Dans tout ce qui suit, nous disons qu'un groupe (fini ou infini) est un {\it 2-groupe élémentaire} si tous ses éléments non triviaux sont d'ordre 2. Bien sûr, les notions de 2-groupe élémentaire et de $\F_2$-espace vectoriel (ainsi que les notions de morphisme correspondantes) coïncident et on ne les distinguera donc pas. On pourra par exemple parler de {\it base} d'un 2-groupe élémentaire, sans faire explicitement référence à sa structure de $\F_2$-espace vectoriel.\\

Soit $E$ un ensemble. On note $Perm(E)$ le groupe des permutations de $E$, ${\cal P}(E)$ l'ensemble des parties de $E$, $[E]$ l'ensemble des parties finies de $E$, et ${\cal E}(E)$ l'ensemble des singletons de $E$. Muni de l'opération de différence symétrique $X + Y = (X\cup Y)\setminus (X\cap Y)$, ${\cal P}(E)$ est un 2-groupe élémentaire (d'élément neutre $\emptyset$) et $[E]$ en est le sous-espace vectoriel de base ${\cal E}(E)$.

\subsection{Généralités sur les groupes de Coxeter.}\label{Généralités sur les groupes de Coxeter.}

Soit $(W,S)$ un système de Coxeter de type $\G = (m_{s,t})_{s,t\in S}$.\\

On note $M_S$ l'ensemble des mots sur $S$ (i.e. le monoïde libre engendré par $S$) et, pour $w\in W$, $M_S(w)$ l'ensemble des éléments de $M_S$ représentant $w$ dans $W$. Les éléments de longueur minimale de $M_S(w)$ sont dit {\it réduits} et leur longueur commune, notée $l_S(w)$, est appelée la {\it longueur} de $w$ (par rapport à $S$). On sait (cf. [B], Ch. IV, \S  1, n° 8, Proposition 7) que l'ensemble des éléments de $S$ apparaissant dans un mot réduit de $M_S(w)$ ne dépend que de $w$ ; cet ensemble est appelé le {\it $S$-support de $w$} et est noté $Supp_S(w)$.\\

Soit $X\subseteq S$. On note $W_X$ le sous-groupe de $W$ engendré par $X$. Alors $W_X$ est constitué des éléments de $W$ dont le $S$-support est inclus dans $X$ (cf. [B], Ch. IV, \S  1, n° 8, Corollaire 1). On en déduit que le couple $(W_X,X)$ est un système de Coxeter de type $(m_{s,t})_{s,t\in X}$, que, pour $X,~Y\subseteq S$, on a $W_X\subseteq W_Y \Longleftrightarrow X\subseteq Y$, et que, si $(X_i)_{i\in I}$ est une famille de parties de $S$, alors $W_{\bigcap_{i\in I}X_i} = \bigcap_{i\in I} W_{X_i}$ (cf. [B], Ch. IV, \S  1, n° 8, Théorème 2). 

Les sous-groupes $W_X$, pour $X\subseteq S$, sont appelés sous-groupes {\it $S$-paraboliques standard} et on dit que $X$ et $W_X$ sont {\it $S$-sphériques} lorsque $W_X$  est fini. Voici un résultat classique (cf. [T], proposition 1, ou [B], exercice 2d page 130) : \\

\begin{propo}\label{tits2}Soit $(W,S)$ un système de Coxeter. Tout sous-groupe fini de $W$ est, à conjugaison près, inclus dans un sous-groupe $S$-sphérique. \\
\end{propo}

Nous notons ${\cal P}_c(S)$ (resp. $[S]_c$) l'ensemble des parties commutatives de $S$ (resp. des parties commutatives et finies de $S$). Les éléments de $[S]_c$ sont des parties $S$-sphériques et, plus généralement, pour $X \in {\cal P}_c(S)$, le sous-groupe $S$-parabolique standard $W_X$ est un 2-groupe élémentaire de base $X$. \\

On note $Aut(\G)$ le groupe des automorphismes de la matrice $\G$ (i.e. le sous-groupe de $Perm(S)$ constitué des permutations de $S$ qui respectent les coefficients $m_{s,t}$). \\

\begin{rema}\label{Rq tilde f}{\bf a.} Tout élément de $Aut(\G)$ se prolonge en un automorphisme de $W$ qui stabilise $S$ et le groupe $Aut(\G)$ s'identifie ainsi au sous-groupe $Aut(W,S)$ de $Aut(W)$.

{\bf b.} Tout élément $f$ de $Perm(S)$ définit, via $X \mapsto f(X) = \{f(x) \mid x\in X\}$, un élément de $Aut([S])$ qui respecte le cardinal et détermine entièrement $f$. Si, de plus, $f$ respecte les coefficients de $\G$, alors $f$ stabilise $[S]_c$ et définit donc un élément de $Aut([S],[S]_c)$ (qui respecte le cardinal). Le groupe $Aut(\G)$ s'identifie ainsi à un sous-groupe de $Aut([S],[S]_c)$.
\end{rema}

\subsection{Le cas à angles droits.}\label{Le cas à angles droits.}

Soit $(W,S)$ un système de Coxeter à angles droits de type $\G = (m_{s,t})_{s,t\in S}$.\\

Dans ce cas, l'ensemble $[S]_c$ des parties commutatives finies de $S$ est exactement l'ensemble des parties $S$-sphériques. Les sous-groupes $S$-sphériques sont donc des 2-groupes élémentaires et, d'après la proposition \ref{tits2}, tous les sous-groupes finis de $W$ sont alors des 2-groupes élémentaires. L'ensemble $F(W)$ est donc constitué de l'élément neutre et des éléments d'ordre 2 de $W$. On voit alors que, pour tout $S' \in {\cal S}(W)$, le système de Coxeter $(W,S')$ est nécessairement à angles droits. Autrement dit, tous les systèmes de Coxeter admis par $W$ sont à angles droits.\\

La proposition suivante, qui apparaît dans la démonstration du corollaire 1 de [T], montre alors en particulier que tous les élements de ${\cal S}(W)$ ont même cardinal, à savoir la dimension du 2-groupe élémentaire $W^{ab}$ sur $\F_2$. \\

Les deux affirmations qui apparaissent dans la définition \ref{def angles droits} sont ainsi prouvées. \\

\begin{propo}\label{isom de [S] sur[S'] envoyant [S]c sur[S']c}Soit $(W,S)$ un système de Coxeter à angles droits. 
\begin{enumerate}
\item L'application ${\cal E}(S)\rightarrow W^{ab},\ \{s\} \mapsto \overline s$ se prolonge en un isomorphisme $\pi_S$ de $[S]$ sur $W^{ab}$. 
\item De plus, $\pi_S$ envoie $[S]_c$ sur $\overline F$.
\end{enumerate}
\end{propo}
\proof D'après la propriété universelle du système de Coxeter $(W,S)$, l'application $S \rightarrow [S]$, $s\mapsto \{s\}$ se prolonge (de façon unique) en un morphisme de groupes de $W$ dans $[S]$. Comme $[S]$ est abélien, ce morphisme passe au quotient en le morphisme $\rho_S : W^{ab}\rightarrow [S]$, donné par $\overline s \mapsto \{s\}$. Comme $\rho_S$ envoie la famille génératrice $(\overline s)_{s\in S}$ du $\F_2$-espace vectoriel $W^{ab}$ sur la base $(\{s\})_{s\in S}$ du $\F_2$-espace vectoriel $[S]$, $\rho_S$ est un isomorphisme, d'où le premier point, avec $\pi_S = (\rho_S)^{-1}$. 

Montrons que $\pi_S([S]_c) = \overline F$. Soit $X\in [S]_c$. Alors $\pi_S(X) = \prod_{s\in X}\overline s$ appartient à $\overline F$, puisque l'élément $\prod_{s\in X}s$ de $W$ est d'ordre 1 ou 2 (donc appartient à $F(W)$). On a donc  $\pi_S([S]_c) \subseteq \overline F$. Réciproquement, soit $w \in F(W)$. D'après la proposition \ref{tits2}, il existe une partie $S$-sphérique $X$ (de $S$) telle que $w$ est conjugué à un élément de $W_X$. Comme $(W,S)$ est à angles droits, on a $X \in [S]_c$ et $w$ est donc conjugué à un élément $\prod_{s\in Y}s$, pour un certain $Y \subseteq X$. On a alors $Y \in [S]_c$ et $\overline w = \prod_{s\in Y}\overline s = \pi_S(Y)$. On a donc $\overline F \subseteq \pi_S([S]_c)$ et le résultat. \qed \\

\begin{rema}\label{relations et singletons}{\bf a.} Si $S'$ est un autre ensemble de Coxeter pour $W$, alors l'application $(\pi_{S'})^{-1}\circ \pi_S$ est un isomorphisme de $[S]$ sur $[S']$ envoyant $[S]_c$ sur $[S']_c$.

{\bf b.} L'isomorphisme $\pi_S$ nous permet d'identifier les groupes $Aut([S])$ et $Aut(W^{ab})$ d'une part, et les groupes $Aut([S],[S]_c)$ et $Aut(W^{ab},\overline F)$ d'autre part. 

{\bf c.} Revenons à la remarque \ref{Rq tilde f} b. On supose ici que $\G$ est à angles droits  ; on montre alors facilement que le groupe $Aut(\G)$ s'identifie précisément au sous-groupe de $Aut([S],[S]_c)$ constitué des éléments de $Aut([S],[S]_c)$ qui respectent le cardinal.\\
\end{rema}

\begin{comms}\label{commentaire Aut([S],[S]c)}Identifions $Aut([S],[S]_c)$ et $Aut(W^{ab},\overline F)$ grâce à $\pi_S$.

{\bf a.} Pour montrer que la suite $\{1\}\hookrightarrow \ker(\pi_{Aut}) \hookrightarrow Aut(W)\stackrel{\pi_{Aut}}{\longrightarrow} Aut(W^{ab},\overline F) \rightarrow \{1\}$ est exacte et scindée, J. Tits a défini la section de $\pi_{Aut}$ qui à tout $\p \in Aut(W^{ab},\overline F) = Aut([S],[S]_c)$ associe l'automorphisme de $W$ donné sur $S$ par $s\mapsto \prod_{x \in \p(\{s\})}x$. 

{\bf b.} Via cette section, $Aut([S],[S]_c)$ s'identifie à un sous-groupe de $Aut(W)$ et on a alors $Aut(W) = \ker(\pi_{Aut})\rtimes Aut([S],[S]_c)$. On vérifie aisément que $Aut([S],[S]_c)$ s'identifie précisément au sous-groupe de $Aut(W)$ constitué des automorphismes de $W$ qui stabilisent l'ensemble des éléments de support $S$-sphérique (i.e. la réunion des sous-groupes $S$-sphériques de $W$).

{\bf c.} Via cette section, on vérifie facilement que le sous-groupe de $Aut([S],[S]_c)$ identifié à $Aut(\G)$ comme en remarque \ref{Rq tilde f} b (ou \ref{relations et singletons} c) s'envoie sur le sous-groupe $Aut(W,S)$ de $Aut(W)$ ; on retrouve ainsi l'identification de la remarque \ref{Rq tilde f} a.
\end{comms}

\section{Outils.}\label{partie outils}

\subsection{Ensembles à relations.}\label{defs}

\begin{defi}[Relations]\label{def relations}Nous appelons {\it ensemble à relation}, ou plus simplement {\it relation}, tout couple $(E,R)$ où $E$ est un ensemble et $R$ une relation binaire symétrique et réflexive sur $E$. Autrement dit, $R$ est une partie de $E\times E$ telle que, pour tous $x,~y \in E$, $(x,x)\in R$, et $(x,y)\in R \Longleftrightarrow (y,x)\in R$. On note $xRy$ pour $(x,y) \in R$.

Soient $(E,R)$ et $(E',R')$ deux ensembles à relation. Un {\it isomorphisme de relations} de $(E,R)$ sur $(E',R')$ est une bijection $f : E\stackrel{\sim}{\longrightarrow} E'$ telle que $(f\times f)(R) = R'$. S'il en existe un, on dit que les relations $(E,R)$ et $(E',R')$ sont isomorphes. En particulier, si $(E,R) = (E',R')$, on parle d'{\it automorphismes} de $(E,R)$ et on note $Aut(E,R)$ le groupe qu'ils constituent.\\
\end{defi}

\begin{nota}[La relation associée à {\boldmath $\G$}]\label{le relation S}Soit $(W,S)$ un système de Coxeter de type $\G = (m_{s,t})_{s,t\in S}$. On munit $S$ de la relation binaire symétrique et réflexive $R_\G = \{(s,t) \in S\times S \mid m_{s,t} = 1 \textrm{ ou } 2\} = \{(s,t) \in S\times S \mid st = ts\}$. \\
\end{nota}

\begin{rema}\label{relation S, cas à angles droits}Soit $W$ un groupe de Coxeter à angles droits.

{\bf a.} Pour $S\in {\cal S}(W)$ et $\G = \G_{(W,S)}$, la relation $R_\G$ détermine entièrement $\G$ puisque, pour $s,~t \in S$ distincts, on a les équivalences suivantes : $m_{s,t} = 2 \Leftrightarrow sR_\G t$, et $m_{s,t} = \infty \Leftrightarrow (s,t) \not \in R_\G$. 

{\bf b.} On voit alors que, pour $S,~S' \in {\cal S}(W)$ et $\G = \G_{(W,S)}$, $\G' = \G_{(W,S')}$, les notions d'isomorphisme (de matrices) de $\G$ sur $\G'$ et d'isomorphisme de relations de $(S,R_\G)$ sur $(S',R_{\G'})$ coïncident. En particulier, on a $Aut(\G) = Aut(S,R_\G)$.\\
\end{rema}

\begin{defi}[Commutants. Parties commutatives]\label{defs parties commutatives}Soit $(E,R)$ un ensemble à relation. Considérons l'application $\fonc{C = C_{(E,R)}}{{\cal P}(E)}{{\cal P}(E)}{X}{\{y\in E \mid \forall x\in X, \ yRx \}}$. 

Pour simplifier, nous notons $C(x)$ l'image $C(\{x\})$ du singleton $\{x\}$ de $E$. Par analogie avec ce qui se passe pour la relation $(S,R_\G)$ associée à une matrice de Coxeter $\G$, nous disons que :
\begin{itemize}
\item pour $x,~y\in E$, $x$ et $y$ {\it commutent} lorsque $xRy$ (i.e $x \in C(y)$), 
\item pour $X \subseteq E$, l'image $C(X)$ est le {\it commutant de $X$ (dans $E$)} et $X$ est {\it commutative} si les éléments de $X$ commutent deux à deux (c'est-à-dire si $X \subseteq C(X)$),
\end{itemize}
et nous notons ${\cal P}_c(E)$ (resp. $[E]_c$) l'ensemble des parties commutatives de $E$ (resp. des parties commutatives et finies de $E$). \\
\end{defi}

\begin{rema}\label{propriétés de C}{\bf a.} $C$ est décroissante (et donc $C^2$ est croissante) pour l'inclusion.

{\bf b.} Pour tout $X \in {\cal P}(E)$, $X\subseteq C^2(X)$.

{\bf c.} On déduit des assertions a et b que $C^3 = C$. 

{\bf d.} Pour tout $X \in {\cal P}(E)$, $X$ est commutative $\Longleftrightarrow C^2(X)$ est commutative.

{\bf e.} Pour tout $X \in {\cal P}(E)$, $C(X) = \bigcap_{x\in X}C(x)$.\\
\end{rema}

\begin{nota}\label{relation cal R}Soit $(E,R)$ un ensemble à relation.

Nous notons $R_{\cal P}$ la relation binaire sur ${\cal P}(E)$ donnée par $XR_{\cal P}Y \Leftrightarrow \forall (x,y) \in X\times Y, \ xRy$. 

On vérifie que $R_{\cal P}$ est symétrique (car $R$ l'est) et que, pour $X \in {\cal P}(E)$, on a $XR_{\cal P}X$ si et seulement si $X$ est commutative.\\
\end{nota}

\begin{rema}\label{rems rel Rp}Conservons les notations introduites ci-dessus. 

{\bf a.} Soit ${\cal X} \subseteq {\cal P}_c(E)$. On note encore $R_{\cal P}$ la relation $R_{\cal P}\cap {\cal X}^2$ induite par $R_{\cal P}$ sur $\cal X$. C'est une relation symétrique et réflexive, et le couple $({\cal X},R_{\cal P})$ est donc un ensemble à relation. 

{\bf b.} Pour $X,~Y \in {\cal P}_c(E)$, on a $XR_{\cal P}Y \Leftrightarrow (X\cup Y) \in {\cal P}_c(E) \Leftrightarrow (X+Y) \in {\cal P}_c(E)$.\\
\end{rema}

Dans ce qui suit, nous serons amenés à considérer différents ensembles de parties commutatives d'un ensemble à relation $(E,R)$ (cf. définitions \ref{defs bicommutants cellules} et \ref{def noyaux}). Cependant nous noterons $R_{\cal P}$, sans plus de précision, les différentes relations associées à ces différents ensembles de parties commutatives ; cela ne provoquera pas de difficulté majeure puisque la relation $XR_{\cal P}Y$ signifiera, dans chaque cas, que la réunion $X\cup Y$ est encore commutative.

\subsection{Cellules et noyaux.}\label{Outils. Premières propriétés.}

Soit $(E,R)$ un ensemble à relation et soit $C = C_{(E,R)}$ (cf. définition \ref{defs parties commutatives}).\\

\begin{defi}[Cellules]\label{defs bicommutants cellules}Nous appelons {\it cellules} (de $(E,R)$) les images par $C^2$ des singletons de $E$ et nous notons ${\cal C}(E,R)$ l'ensemble qu'elles constituent. Nous disons plus précisément que la cellule $C^2(x)$ est la {\it cellule de $x$} et que $x$ {\it définit} la cellule $C^2(x)$. \\
\end{defi}

\begin{comms}\label{cellules Tits}
Soit $X \in {\cal P}_c(E)$. On montre que $C^2(X)$ est l'intersection des parties commutatives maximales de $E$ contenant $X$. En particulier, pour $x \in E$, la cellule $C^2(x)$ est l'intersection des parties commutatives maximales de $E$ contenant $x$. 
C'est sous cette deuxième forme que la notion de  cellule d'un élément a été introduite par J. Tits dans la partie finale de [T] ($C^2(x)$ y est noté $T(x)$). Lorsque $S$ est infini, la définition que l'on donne ici a l'avantage d'éviter d'avoir à utiliser le lemme de Zorn (pour l'existence des parties commutatives maximales). \\
\end{comms}

\begin{rema}\label{remarque bicommutant}{\bf a.} La remarque \ref{propriétés de C} b nous montre que, pour tout $x\in E$, on a $x \in C^2(x)$, et la remarque \ref{propriétés de C} d nous montre que les cellules sont des parties commutatives de $E$.

{\bf b.} Comme $C^3 = C$ (remarque \ref{propriétés de C} c), les cellules sont des points fixes de $C^2$ et on a donc en particulier, pour toute cellule $T$, $x\in T \Longrightarrow C^2(x) \subseteq T$.\\
\end{rema}

\begin{defi}[Noyaux]\label{def noyaux}Notons $\equiv$ la relation d'équivalence sur $E$ donnée par $y\equiv z \Longleftrightarrow C^2(y) = C^2(z)$ (ou encore $y\equiv z \Longleftrightarrow C(y) = C(z)$, puisque $C^3 = C$). Nous appelons {\it noyaux} (de $(E,R)$) les classes d'équivalence de $E$ pour $\equiv$ et nous notons ${\cal N}(E,R)$ l'ensemble qu'ils constituent.\\
\end{defi}

\begin{nota}Soit $T\in {\cal C}(E,R)$. Il est clair que l'ensemble $\{x\in E \mid C^2(x) = T\}$ est un noyau de $(E,R)$ ; on l'appelle {\it le noyau de $T$} et on le note $Noy(T)$ (autrement dit, $Noy(T) = N(x)$ pour tout $x \in E$ tel que $T = C^2(x)$).

Soit $N \in {\cal N}(E,R)$. On appelle {\it cellule de $N$}, et on note $Cel(N)$, la cellule commune aux éléments de $N$ (autrement dit, $Cel(N) = C^2(x)$ pour tout $x\in N$). \\
\end{nota}

\begin{rema}\label{props equiv}{\bf a.} Soit $N$ un noyau. On a $N \subseteq Cel(N)$ (grâce à la remarque \ref{remarque bicommutant} a). En particulier, les noyaux sont des parties commutatives (i.e. ${\cal N}(E,R) \subseteq {\cal P}_c(E)$). 

{\bf b.} Toute cellule est la réunion (disjointe) des noyaux qu'elle rencontre. 

{\bf c.} Les applications $Cel : {\cal N}(E,R) \rightarrow {\cal C}(E,R)$ et $Noy : {\cal C}(E,R) \rightarrow {\cal N}(E,R)$ sont des bijections inverses l'une de l'autre.\\
\end{rema}

Notons que, les cellules et les noyaux étant des parties commutatives de $E$, on peut munir chacun des ensembles ${\cal C}(E,R)$ et ${\cal N}(E,R)$ d'une structure d'ensemble à relation, avec la relation ${R_{\cal P}}$ définie en notation \ref{relation cal R} (cf. remarque \ref{rems rel Rp} a). Nous avons alors la proposition suivante : \\

\begin{propo}\label{Cellules et noyaux d'un relation régulier}Les bijections $Cel$ et $Noy$ décrites ci-dessus sont des isomorphismes de relations inverses l'un de l'autre.
\end{propo}
\proof Il s'agit de voir que, pour $N, P \in {\cal N}(E,R)$, on a $N{R_{\cal P}}P \Leftrightarrow Cel(N){R_{\cal P}}Cel(P)$. Comme $N\subseteq  Cel(N)$ et $P\subseteq  Cel(P)$, l'implication $Cel(N){R_{\cal P}} Cel(P) \Rightarrow N{R_{\cal P}}P$ est évidente. 

Réciproquement, supposons $N {R_{\cal P}} P$. Posons $Cel(N) = C^2(x)$ et $Cel(P) = C^2(y)$, pour $x \in N$ et $y \in P$. On a $xRy$, c'est-à-dire $x \in C(y)$, donc $C^2(y) \subseteq C(x)$, et on veut montrer que, pour tout $x' \in Cel(N) = C^2(x)$ et tout $y'\in Cel(P) = C^2(y)$, on a $x'Ry'$. Comme $y'\in C^2(y)\subseteq C(x)$, on a $C^2(x)\subseteq C(y')$, et comme $x' \in C^2(x)$, on a $x' \in C(y')$, c'est-à-dire $x'Ry'$, d'où le résultat. \qed\\

Donnons à présent une propriété importante des noyaux, qui nous servira à la fois pour montrer qu'un groupe de Coxeter à angles droits $W$ est rigide (cf. section \ref{sur les isomorphismes ...}, théorème \ref{theorème principal}) et pour étudier le groupe $Aut([S],[S]_c)$ (cf. section \ref{premier dévissage}, proposition \ref{premières décompositions}).\\

\begin{lem}\label{prop noyaux}Pour $N, \ P \in {\cal N}(E,R)$, on a : $ N{R_{\cal P}} P \Longleftrightarrow \exists\,  (x,y)\in N\times P,\  xRy$.
\end{lem}
\proof Le sens direct étant évident (par définition de $R_{\cal P}$), montrons la réciproque.

 Soient $x \in N$, $y \in P$ tels que $xRy$, et soient $x' \in N$, $y' \in P$. On a $C(x) = C(x')$ et $C(y) = C(y')$ par définition de $N$ et $P$, et $x \in C(y)$ par hypothèse sur $x$ et $y$. On a alors $x \in C(y')$, ou encore $y' \in C(x)$ (puisque $R$ est symétrique), d'où $y' \in C(x')$, c'est-à-dire $x'Ry'$. On a donc le résultat.\qed\\

\begin{propo}\label{isom entre noyaux}Soient $(E,R)$ et $(E',R')$ deux ensembles à relation. Supposons qu'il existe un isomorphisme de relations $\psi$ respectant le cardinal, de $({\cal N}(E,R),{R_{\cal P}})$ sur $({\cal N}(E',R'),{R'_{\cal P}})$. 

Alors toute bijection $f$ de $E$ sur $E'$ définie, noyau par noyau, par n'importe quelle bijection de $N \in {\cal N}(E,R)$ sur $\psi (N) \in {\cal N}(E',R')$ est un isomorphisme de relations de $(E,R)$ sur $(E',R')$. 

En particulier, les ensembles à relation $(E,R)$ et $(E',R')$ sont isomorphes.
\end{propo}
\proof  Remarquons que la construction proposée est possible (grâce à l'axiome du choix) car ${\cal N}(E,R)$ (resp. ${\cal N}(E',R')$) est une partition de $E$ (resp. $E'$) et car les noyaux $N$ et $\psi(N)$ ont même cardinal par hypothèse.

Soient $x, ~y \in E$. Notons $N$ (resp. $P$) le noyau de $x$ (resp. $y$). Grâce au lemme \ref{prop noyaux} et au fait que $\psi$ est un isomorphisme de relations, on a les équivalences suivantes : $xRy \Longleftrightarrow N {R_{\cal P}}P \Longleftrightarrow \psi(N)R'_{\cal P}\psi(P) \Longleftrightarrow f(x)R'f(y)$. La bijection $f$ est donc un isomorphisme de relations.\qed 

\subsection{Théorème principal.}\label{sur les isomorphismes ...}

Soient $(E,R)$ et $(E',R')$ deux ensembles à relations. Posons, pour simplifier, $C = C_{(E,R)}$, ${\cal C} = {\cal C}(E,R)$, ${\cal N} = {\cal N}(E,R)$ et $C' = C_{(E',R')}$, ${\cal C}' = {\cal C}(E',R')$, ${\cal N}' = {\cal N}(E',R')$.\\

Supposons qu'il existe un isomorphisme $\F_2$-linéaire $\p$ de $[E]$ sur $[E']$ envoyant $[E]_c$ sur $[E']_c$.\\

\begin{lem}\label{phi isom de rels}Soient $X,~Y \in [E]_c$. Alors $\p(X),~\p(Y)\in [E']_c$ et on a : 
\begin{center}$X\cup Y$ est commutative $\Longleftrightarrow \p(X) \cup \p(Y)$ est commutative.\end{center}
\end{lem}
\proof On suppose $X,~Y \in [E]_c$ et $\p([E]_c) = [E']_c$, donc $\p(X),~\p(Y)\in [E']_c$. De plus, d'après la remarque \ref{rems rel Rp} b, on a $X\cup Y \in [E]_c \Leftrightarrow X+Y \in [E]_c$, et $\p(X) \cup \p(Y)\in [E']_c \Leftrightarrow \p(X) + \p(Y) \in [E']_c$. Enfin, comme $\p$ est un isomorphisme linéaire envoyant $[E]_c$ sur $[E']_c$, on a $X+Y \in [E]_c \Leftrightarrow \p(X)+\p(Y) \in [E']_c$ et on en déduit le résultat.  \qed \\

\begin{nota}Si $T$ est une cellule d'un ensemble à relation $(E,R)$, nous posons $T^\star = T\setminus Noy(T) = \{x\in E \mid C^2(x) \varsubsetneq T\}$ et nous notons $V_T$ l'espace vectoriel quotient $[T]/[T^\star]$.\\
\end{nota}

\begin{propo}\label{sous-espaces [T] permutés}Conservons les notations introduites ci-dessus. 
\begin{enumerate}
\item Soit $x\in E$. Il existe $x'\in \p(\{x\})$ tel que $x \in \p^{-1}(\{x'\})$ et, pour tout tel $x'$, on a $\p([C^2(x)]) = [C'^2(x')]$.
\item Soit $T \in {\cal C}$ et soit $T'$ l'unique élément de ${\cal C}'$ tel que $\p([T]) = [T']$. Alors $\p([T^\star]) = [T'^\star]$ et $\p$ induit donc un isomorphisme $\bar \p_T : V_T \stackrel{\sim}{\longrightarrow} V_{T'}$ donné par $X + [T^\star] \longmapsto \p(X) + [T'^\star]$.
\end{enumerate}
\end{propo}
\proof Montrons la première assertion. L'existence de $x'\in \p(\{x\})$ tel que $x \in \p^{-1}(\{x'\})$ résulte des égalités $\{x\} = \p^{-1}(\p(\{x\})) = \p^{-1}(\sum_{x'\in \p(\{x\})}\{x'\}) = \sum_{x'\in \p(\{x\})}\p^{-1}(\{x'\})$ et du fait que l'on a $\Sigma_{x'\in \p(\{x\})}\p^{-1}(\{x'\})\subseteq \bigcup_{x'\in \p(\{x\})}\p^{-1}(\{x'\})$.

Pour montrer l'inclusion $\p([C^2(x)]) \subseteq  [C'^2(x')]$, il suffit de montrer que, pour tout $y\in C^2(x)$, on a $\p(\{y\}) \in [C'^2(x')]$, c'est-à-dire que, pour tout $y' \in C'(x')$, la partie $\{y'\} \cup \p(\{y\})$ est commutative. Soient donc $y\in C^2(x)$ et $y' \in C'(x')$. Les éléments $x'$ et $y'$ commutent (i.e. $\{x'\}\cup \{y'\}$ est commutative), donc, par le lemme précédent, $\p^{-1}(\{x'\})\cup \p^{-1}(\{y'\})$ est une partie commutative, qui de plus contient $x$ (par hypothèse sur $x'$). On en déduit que l'on a $\p^{-1}(\{y'\}) \subseteq C(x)$, donc que $\p^{-1}(\{y'\})\cup \{y\}$ est commutative (puisque $y \in C^2(x)$), et donc (grâce au lemme précédent) que $\{y'\} \cup \p(\{y\})$ est également commutative. Puisque les rôles joués par $(\p,x)$ et $(\p^{-1},x')$ sont symétriques, on obtient de la même façon l'inclusion $\p^{-1}([C'^2(x')]) \subseteq  [C^2(x)]$ et donc $\p([C^2(x)]) = [C'^2(x')]$.  

Montrons la seconde assertion. D'après le premier point, pour $T \in {\cal C}$, il existe $T' \in {\cal C}'$ tel que $\p([T]) = [T']$, et un tel $T'$ est nécessairement unique puisque si $X_1$ et $X_2$ sont deux ensembles tels que $[X_1] = [X_2]$, alors on a clairement $X_1 = X_2$.

Soit $y \in T^\star$. D'après le premier point, $\p([C^2(y)])$ est l'ensemble des parties finies d'une cellule $U'$ de $(E',R')$ nécessairement strictement incluse dans $T'$, puisque $C^2(y) \varsubsetneq T$ (donc $[C^2(y)] \varsubsetneq [T]$). On a donc $\p(\{y\}) \in [U'] \subseteq [T'^\star]$, ce qui montre l'inclusion $\p([T^\star]) \subseteq [T'^\star]$. L'autre inclusion se démontre de la même façon, à partir de $\p^{-1}$, et on a donc $\p([T^\star]) = [T'^\star]$. La fin de l'assertion s'en déduit immédiatement. \qed \\

\begin{nota}\label{nota phiC}On déduit de la proposition précédente que l'isomorphisme linéaire $\p : [E] \stackrel{\sim}{\longrightarrow} [E']$ envoie de façon bijective (en respectant l'inclusion et la dimension) l'ensemble des sous-espaces de $[E]$ de la forme $[T]$, où $T\in {\cal C}$, sur l'ensemble des sous-espaces de $[E']$ de la forme $[T']$, où $T'\in {\cal C}'$. L'isomorphisme $\p$ induit donc une bijection $\p_{\cal C}$ de ${\cal C}$ sur ${\cal C}'$ donnée par : 
$$\bij{\p_{\cal C}}{\cal C}{{\cal C}'}{T}{T'},\ \textrm{  où   } \p([T]) = [T'].$$
\end{nota}

\begin{propo}\label{prop phiC}Conservons les notations introduites ci-dessus. 

La bijection $\p_{\cal C}$ est un isomorphisme de relations de $({\cal C},R_{\cal P})$ sur $({\cal C}',R'_{\cal P})$.
\end{propo}
\proof Soient $T,~U \in {\cal C}$ et soient $T' = \p_{\cal C}(T)$, $U' = \p_{\cal C}(U) \in {\cal C}'$. Il s'agit de montrer que l'on a $TR_{\cal P}U \Leftrightarrow T'R'_{\cal P}U'$, c'est-à-dire que $T\cup U$ est commutative si et seulement si $T'\cup U'$ l'est. Supposons $T'\cup U'$ commutative et soit $(t,u)\in T\times U$. La partie $\p(\{t\}) \cup \p(\{u\})$ de $E'$ est incluse dans $T'\cup U'$ donc est commutative et, d'après le lemme \ref{phi isom de rels}, on a donc $tRu$. Comme ceci est vrai pour tout $(t,u)\in T\times U$, on a $TR_{\cal P}U$ et l'implication $T'R'_{\cal P}U' \Rightarrow TR_{\cal P}U$ est donc démontrée. L'autre implication se démontre de la même façon, en utilisant $\p^{-1}$.\qed \\

Donnons à présent le point clé de la démonstration du théorème \ref{Les groupes de Coxeter à angles droits sont rigides} : \\

\begin{theo}\label{theorème principal}Soient $(E,R)$ et $(E',R')$ deux ensembles à relations et soit $\p$ un isomorphisme de $[E]$ sur $[E']$ envoyant $[E]_c$ sur $[E']_c$.
\begin{enumerate}
\item L'isomorphisme $\p$ induit un isomorphisme de relations qui respecte le cardinal $$\fonc{\p_{\cal N}}{({\cal N},{R_{\cal P}}) }{({\cal N}',{R'_{\cal P}})}{N}{N'},\ \textrm{  où  } \p([Cel(N)]) = [Cel(N')].$$

\item Les ensembles à relation $(E,R)$  et $(E',R')$ sont isomorphes.
\end{enumerate}
\end{theo}
\proof 
Le fait que l'application $\p_{\cal N}$ définie en 1 soit un isomorphisme de relations résulte des propositions \ref{Cellules et noyaux d'un relation régulier} et \ref{prop phiC}, puisque l'on a $\p_{\cal N} = Noy \circ \p_{\cal C} \circ Cel$ avec les isomorphismes de relations $Cel : ({\cal N},{R_{\cal P}}) \stackrel{\sim}{\longrightarrow} ({\cal C},{R_{\cal P}})$, $\p_{\cal C} : ({\cal C},{R_{\cal P}}) \stackrel{\sim}{\longrightarrow} ({\cal C}',{R'_{\cal P}})$ et $Noy : ({\cal C}',{R'_{\cal P}})\stackrel{\sim}{\longrightarrow} ({\cal N}',{R'_{\cal P}})$.

Montrons que $\p_{\cal N}$ respecte le cardinal. Soient $N \in {\cal N}$ et $N' = \p_{\cal N}(N)\in {\cal N}'$. Posons $T = Cel(N)$ et $T' = Cel(N')$. On a $\p([T] = [T']$. D'après la seconde assertion de la proposition \ref{sous-espaces [T] permutés}, les espaces vectoriels $V_{T}$ et $V_{T'}$ sont isomorphes (via $\bar \p_{T}$) et ont donc même dimension. Or $T$ est la réunion disjointe de $N$ et de $T^\star$, donc $[T] = [N]\oplus [T^\star]$ et les espaces vectoriels $[N]$ et $V_{T} = [T]/[T^\star]$ sont donc isomorphes. De même, les espaces vectoriels $[N']$ et $V_{T'}$ sont isomorphes et on a le résultat puisque l'espace vectoriel $[N]$ (resp. $[N']$), qui admet pour base l'ensemble des singletons de $N$ (resp. de $N'$) est de dimension $Card(N)$ (resp. $Card(N')$).

Le second point résulte du premier et de la proposition \ref{isom entre noyaux}.\qed 

\section{Sur la rigidité des groupes de Coxeter à angles droits.}\label{Problèmes d'isomorphie}

\subsection{Les groupes de Coxeter à angles droits sont rigides.}\label{Conclusion rigidité.}

\begin{defi}[Rigidité]\label{defs rigidité}Soit $W$ un groupe de Coxeter. 

On dit que $W$ est {\it rigide} lorsque tous les types admis par $W$ sont isomorphes, autrement dit lorsque, pour tous $S,~S' \in {\cal S}(W)$, il existe $\a \in Aut(W)$ tel que $\a(S) = S'$.\\
\end{defi}

Les résultats obtenus dans la partie précédente nous permettent de démontrer le théorème \ref{Les groupes de Coxeter à angles droits sont rigides} annoncé en introduction : \\

\begin{theo}\label{Les groupes de Coxeter à angles droits sont rigides}Les groupes de Coxeter à angles droits sont rigides.
\end{theo}
\proof Soit $W$ un groupe de Coxeter à angles droits et soient $S, ~ S'\in {\cal S}(W)$. Notons $\G$ et $\G'$ les types respectifs des systèmes $(W,S)$ et $(W,S')$ ; ce sont des matrices de Coxeter à angles droits (cf. section \ref{Le cas à angles droits.}). Considérons les relations $(S,R)$ et $(S',R')$ associées aux systèmes $(W,S)$ et $(W,S')$ respectivement : on a $R = \{(s,t) \in S\times S \mid st = ts\}$ et $R' = \{(s',t') \in S'\times S' \mid s't' = t's'\}$.

D'après la proposition \ref{isom de [S] sur[S'] envoyant [S]c sur[S']c}, il existe un isomorphisme de $[S]$ sur $[S']$ envoyant $[S]_c$ sur $[S']_c$ (par exemple l'isomorphisme $(\pi_{S'})^{-1} \circ \pi_S$ de la remarque \ref{relations et singletons} a). On en déduit, grâce au théorème \ref{theorème principal}, que les relations $(S,R)$ et $(S',R')$ sont isomorphes. Comme les matrices $\G$ et $\G'$ sont à angles droits, un isomorphisme de relations de $(S,R)$ sur $(S',R')$ est un isomorphisme de $\G$ sur $\G'$ (cf. remarque \ref{relation S, cas à angles droits} b). On a donc le résultat. \qed

\subsection{Sur la rigidité forte.}\label{Notions de rigidité.}

\begin{defi}[Rigidité forte]\label{def rigidité forte}Soit $W$ un groupe de Coxeter.

On dit que $W$ est {\it fortement rigide} lorsque, pour tous $S, ~S' \in {\cal S}(W)$, il existe $w \in W$ tel que $wSw^{-1} = S'$ (c'est-à-dire qu'il existe un automorphisme intérieur envoyant $S$ sur $S'$).\\
\end{defi}

\begin{rema}\label{notions de rigidité}Soient $W$ un groupe de Coxeter, $S\in {\cal S}(W)$ et $\G = \G_{(W,S)}$. Identifions $Aut(\G)$ et $Aut(W,S)$ (comme en remarque \ref{Rq tilde f} a). 

{\bf a.} On vérifie facilement que l'on a la caractérisation  : 
$$W  \textrm{  est fortement rigide  } \Longleftrightarrow W  \textrm{  est rigide et  } Aut(W) = Int(W) \cdot Aut(\G).$$

{\bf b.} Supposons $W$ à angles droits.

Alors $W$ est rigide, d'après le théorème \ref{Les groupes de Coxeter à angles droits sont rigides}. De plus, on a $Aut(W) = \ker(\pi_{Aut}) \rtimes Aut([S],[S]_c)$ ([T], Corollaire 1, ou commentaires \ref{commentaire Aut([S],[S]c)} b  et \ref{commentaire Aut([S],[S]c)} c ci-dessus) , avec $Aut(\G) \subseteq Aut([S],[S]_c)$ (remarque \ref{Rq tilde f} b ou \ref{relations et singletons} c) et clairement $Int(W) \subseteq \ker(\pi_{Aut})$. La caractérisation de la rigidité forte de l'assertion a ci-dessus se traduit donc dans ce cas par :
$$W \textrm{  est fortement rigide  } \Longleftrightarrow \ker(\pi_{Aut}) = Int(W) \textrm{  et  } Aut(([S],[S]_c)) = Aut(\G).$$\\
\end{rema}

\begin{comms}\label{commentaires rigidité forte}Soit $(W,S)$ un système de Coxeter à angles droits de type $\G$. 

{\bf a.} Notons ${^W\negthickspace S}$ l'ensemble des conjugués des éléments de $S$ dans $W$. Le corollaire 1 de [T] montre en particulier que l'on a $\ker(\pi_{Aut}) \subseteq Aut(W,{^W\negthickspace S})$. De plus, on vérifie que l'on a $Aut([S],[S]_c)\cap Aut(W,{^W\negthickspace S}) = Aut(\G)$. On en déduit la caractérisation $Aut(([S],[S]_c)) = Aut(\G) \Leftrightarrow Aut(W) = Aut(W,{^W\negthickspace S})$.

{\bf b.} Lorsque $W$ est de rang {\it fini}, on peut, par des considérations simples sur $\G$ (ou sur son graphe), déterminer si $W$ est fortement rigide ou non. En effet, on a :  
\begin{enumerate}
\item $\ker(\pi_{Aut}) = Int(W) \Leftrightarrow \forall s\in S,\ \forall t,~u \in S\setminus C(s), \ \exists \ t_0 = t,~t_1,~\ldots,~t_n = u \in S\setminus C(s)$ tels que $t_{i-1}$ et $t_i$ commutent, pour $1\leq i \leq n$ (cf. [M], corollaire du théorème principal),
\item $Aut(([S],[S]_c)) = Aut(\G) \Leftrightarrow \forall s\in S, \ C^2(s) = \{s\}$ (cf. [BM], théorème 5.1, appliqué aux groupes de Coxeter à angles droits et de rang fini, et les commentaires \ref{cellules Tits} et \ref{commentaires rigidité forte} a).
\end{enumerate}
Ces conditions apparaissent aussi dans [BMMN] (théorème 4.10).

{\bf c.} Lorsque $W$ est à angles droits et de rang infini, la condition énoncée en b 1 est nécessaire, pour avoir $\ker(\pi_{Aut}) = Int(W)$, mais n'est pas suffisante (cf. [T], proposition 5 et remarque finale de la partie 3). \\
\end{comms}

La proposition \ref{caract} ci-dessous montre que la caractérisation du commentaire \ref{commentaires rigidité forte} b 2 est encore valable pour les groupes de Coxeter à angles droits de rang {\it infini}. C'est un résultat qui découle également de l'étude générale du groupe $Aut(([S],[S]_c))$ que nous effectuons dans la partie suivante (cf. remarque \ref{remarque finale} c ci-dessous).\\

\begin{nota}\label{automs élémentaires}Soient $s \in S$ et $t \in C^2(s)$, $t\not= s$. Nous notons $\a_{s,t}$ l'endomorphisme $\F_2$-linéaire de $[S]$ donné par $\a_{s,t}(\{x\}) = \{x\}$ si $x\in S\setminus \{ s\}$, et $\a_{s,t}(\{s\}) = \{s,t\}$. \\
\end{nota}

\begin{rema}\label{rems automs éléms}{\bf a.} Comme $\a_{s,t}$ est clairement involutif, c'est un élément de $Aut([S])$.

{\bf b.} On a en fait $\a_{s,t} \in Aut([S],[S]_c)$. En effet, pour $X \in [S]_c$, on a soit $\a_{s,t}(X) = X \in [S]_c$ si $s \not \in X$, soit $\a_{s,t}(X) = X + \{t\}$ si $s\in X$, auquel cas $t \in C^2(s) \subseteq C^2(X) \subseteq C(X)$ (la dernière inclusion est vérifiée car, $X$ étant commutative, on a $X \subseteq C(X)$) et donc $X + \{t\} \in [S]_c$ ; ceci montre l'inclusion $\a_{s,t}([S]_c) \subseteq [S]_c$, et comme $\a_{s,t}$ est involutif, on a le résultat.\\
\end{rema}

\begin{comms}{\bf a.} La matrice de $\a_{s,t}$ dans la base ${\cal E}(S) = \{\{x\} \mid x \in S\}$ de $[S]$ est une matrice de transvection élémentaire.

{\bf b.} Les automorphismes $\a_{s,t}$, vus comme éléments de $Aut(W)$, apparaissent sous une forme plus générale dans [BM] (lemme 6.1).\\
\end{comms}

\begin{propo}\label{caract}Soit $(W,S)$ un système de Coxeter à angles droits de type $\G = \G_{(W,S)}$. On a $$Aut([S],[S]_c) = Aut(\G) \Longleftrightarrow \forall s\in S, \  C^2(s) = \{s\}.$$
\end{propo}
\proof Supposons que, pour tout $s\in S$, $C^2(s) = \{s\}$. La proposition \ref{sous-espaces [T] permutés} nous montre alors que les sous-espaces de $[S]$ de la forme $[\{s\}] = \{\emptyset,\{s\}\}$, où $s \in S$, sont permutés entre eux par les éléments de $Aut([S],[S]_c)$. Comme ces éléments sont des automorphismes de $[S]$, ils fixent $\emptyset$ (élément neutre de $[S]$) et permutent donc entre eux les singletons de $S$. On en déduit que $Aut([S],[S]_c) \subseteq Aut(\G)$ (cf. remarque \ref{relations et singletons} c).

Réciproquement, supposons qu'une cellule $C^2(s)$ contienne un élément $t \not = s$, et considérons l'endomorphisme $\a_{s,t}$ de $[S]$ défini en notation \ref{automs élémentaires}. C'est un élément de $Aut([S],[S]_c)$, d'après la remarque \ref{rems automs éléms} b, qui n'appartient pas à $Aut(\G)$, puisqu'il ne respecte pas le cardinal. On a donc $Aut(\G) \varsubsetneq Aut([S],[S]_c)$ et le résultat. \qed 

\section{Le groupe {\boldmath $Aut(W^{ab},\overline F) = Aut([S],[S]_c)$}.}\label{description de Aut(F)}

Soit $W$ un groupe de Coxeter à angles droits. Le but de cette partie est d'étudier le groupe $Aut(W^{ab},\overline F)$ intervenant dans la décomposition $Aut(W) = \ker(\pi_{Aut}) \rtimes Aut(W^{ab},\overline F)$ établie par J. Tits dans [T]. \\

Fixons $S \in {\cal S}(W)$ et notons $\G = \G_{(W,S)}$ le type de $(W,S)$. Nous identifions $Aut(W^{ab},\overline F)$ à $Aut([S],[S]_c)$ comme en remarque \ref{relations et singletons} b. C'est précisément ce groupe $Aut([S],[S]_c)$ que nous étudions dans ce qui suit, notamment grâce aux résultats de la section \ref{sur les isomorphismes ...} (appliqués au cas particulier où $S = S'$). \\

Pour alléger les énoncés, nous posons $R = R_\G$, $C = C_{(S,R)}$, ${\cal C} = {\cal C}(S,R)$ et ${\cal N} = {\cal N}(S,R)$.

\subsection{Les sous-groupes {\boldmath $Aut(\G)$} et {\boldmath $K(\G)$}. Dévissage de {\boldmath $Aut([S],[S]_c)$}.}\label{premier dévissage}

Notons $Aut({\cal N},{R_{\cal P}},Card)$ le sous-groupe de $Aut({\cal N},{R_{\cal P}})$ constitué des automorphismes de la relation $({\cal N},{R_{\cal P}})$ qui respectent le cardinal. Le théorème \ref{theorème principal} nous montre en particulier que tout $\p \in Aut([S],[S]_c)$ induit un élément $\p_{\cal N}$ de $Aut({\cal N},{R_{\cal P}},Card)$, donné, pour $N\in {\cal N}$, par $N\mapsto N' (\in {\cal N})$, où $\p([Cel(N)]) = [Cel(N')]$.\\

\begin{nota}[Le groupe {\boldmath $K(\G)$}]Notons $\t$ l'application $ Aut([S],[S]_c)\rightarrow Aut({\cal N},{R_{\cal P}},Card)$, $\p \mapsto \p_{\cal N}$. On vérifie facilement que $\t$ est un morphisme de groupes. Nous notons $K(\G)$ son noyau.\\
\end{nota}

Rappelons que l'on a $Aut(\G) = Aut(S,R)$ (cf. remarque \ref{relation S, cas à angles droits} a). On a vu en remarque \ref{relations et singletons} c que le groupe $Aut(\G)$ s'identifie, via $\s \longmapsto (X \mapsto \s(X))$, au sous-groupe de $Aut([S],[S]_c)$ constitué des éléments de $Aut([S],[S]_c)$ qui respectent le cardinal. On note encore $\s$ l'élement $X \mapsto \s(X)$ de $Aut([S],[S]_c)$.

Soit $\s \in Aut(\G)$. On vérifie que, pour tout $X \in [S]$, on a $\s([X]) = [\s(X)]$ et $C(\s(X)) = \s(C(X))$. L'automorphisme $\s$ permute donc les noyaux (de même cardinal) de $(S,R)$ et l'on voit que l'élément $\t(\s) = \s_{\cal N}$ de $Aut({\cal N},{R_{\cal P}},Card)$ est simplement donné par $N \mapsto \s(N)$.\\

\begin{nota}Soit $\Omega$ l'ensemble des classes d'équivalence de ${\cal N}$ pour la relation $Card(N) = Card(P)$. Pour $\omega \in \Omega$, fixons une fois pour toutes un représentant $N_\omega$ de la classe $\omega$ et, pour tout $N\in \omega$, une bijection $\s_{\omega,N} : N_\omega \longrightarrow N$ (il en existe une, puisque $Card(N_\omega) = Card(N)$).

Pour $\psi \in Aut({\cal N},{R_{\cal P}},Card)$, nous définissons la permutation $\s_\psi$ de $S$, noyau par noyau, de la manière suivante : si $N \in {\cal N}$ et si $\omega$ est la classe de $N$ (et de $\psi(N)$), alors $\s_\psi$ est donné de $N$ sur $\psi(N)$ par $\s_{\omega,\psi(N)}\circ (\s_{\omega,N})^{-1}$. On note $G$ l'ensemble des $\s_\psi$, pour $\psi \in Aut({\cal N},{R_{\cal P}},Card)$.\\
\end{nota}

\begin{propo}\label{premières décompositions}Pour $\psi \in Aut({\cal N},{R_{\cal P}},Card)$, on a $\s_\psi \in Aut(\G)$. De plus, l'application $\psi \mapsto \s_\psi$ est un morphisme de groupes et une section de $\t$. En particulier, $G$ est un sous-groupe de $Aut(\G)$ et on a les décompositions : 
$$Aut([S],[S]_c) = K(\G) \rtimes G\ \textrm{  et  } \ Aut(\G) = (K(\G)\cap Aut(\G)) \rtimes G$$
\end{propo}
\proof La proposition \ref{isom entre noyaux} nous montre que, pour $\psi \in Aut({\cal N},{R_{\cal P}},Card)$, la  permutation $\s_\psi$ de $S$ appartient à $Aut(\G)$. De plus, vu la définition des $\s_\psi$, l'application $Aut({\cal N},{R_{\cal P}},Card) \rightarrow Aut(\G),\ \psi \mapsto \s_\psi$ est un morphisme de groupes et $\s_\psi$ s'envoie sur $\psi$ par $\t$, d'où le résultat. \qed \\


Dans les sections suivantes, nous explicitons le groupe $K(\G)$. Décrivons pour le moment le groupe $K(\G)\cap Aut(\G)$ : \\

\begin{propo}\label{intersection de Aut(G) et K}Soient $\s \in K(\G)\cap Aut(\G)$ et $N \in {\cal N}$. Alors l'automorphisme $\s$ induit (par restriction) une permutation $\s_{| N}$ de $N$. 

De plus, l'application $\rho_1 : K(\G)\cap Aut(\G) \rightarrow \prod_{N\in {\cal N}}Perm(N)$, $\s \mapsto (\s_{| N})_{N\in {\cal N}}$ est un isomorphisme de groupes.
\end{propo}
\proof On sait que tout élément $\s$ de $Aut(\G)$ permute les noyaux et induit donc une bijection de $N$ sur $\s(N)$. Si, de plus, $\s \in K(\G)$, alors $\s(N) = N$ et la permutation $\s_{|N}$ est bien définie. On vérifie facilement que $\s \longmapsto (\s_{| N})_{N\in {\cal N}}$ est un morphisme de groupes injectif (puisque ${\cal N}$ est une partition de $S$). De plus, si $(\s_{N})_{N\in {\cal N}}$ est un élément de $\prod_{N\in {\cal N}}Perm(N)$, alors, d'après la proposition \ref{isom entre noyaux}, la permutation de $S$ définie noyau par noyau par $\s_N : N \longrightarrow N$, pour tout $N\in {\cal N}$, appartient à $Aut(\G)$ (c'est le cas particulier où $\psi = Id_{{\cal N}}$). On voit donc que le morphisme $\s \longmapsto (\s_{| N})_{N\in {\cal N}}$ est surjectif. \qed

\subsection{Les sous-groupes {\boldmath $D(\G)$} et {\boldmath $K^\circ(\G)$}. Dévissage de {\boldmath $K(\G)$}.}\label{le groupe K}

Comme ${\cal N}$ est une partition de $S$, on a $[S] = \oplus_{N\in {\cal N}}[N]$.\\

\begin{propo}\label{caract de F et K}
L'ensemble $[S]_c$ est l'ensemble des sommes (finies) $\sum_{N \in {\cal N}}X_N$ où, pour tout $N \in {\cal N}$, $X_N \in [N]$, et où $\{N \in {\cal N} \mid X_N \not= \emptyset\}$ est fini et de réunion commutative (i.e. formé de noyaux deux à deux en relation ${R_{\cal P}}$). On en déduit que l'on a :
$$K(\G) = \{\p \in Aut([S]) \mid \forall ~T\in {\cal C}, \ \p([T]) = [T]\}.$$
\end{propo}
\proof Toute telle somme $\sum_{N \in {\cal N}}X_N$ est une partie commutative finie de $S$ et appartient donc à $[S]_c$. Réciproquement, si $X \in [S]_c$, alors, dans la décomposition $[S] = \oplus_{N\in {\cal N}}[N]$, $X$ s'écrit $X = \sum_{1\leq k \leq n}X_k$, où, pour tout $k$, $X_k$ est une partie finie non vide d'un noyau $N_k$ et où, d'après le lemme \ref{prop noyaux}, les noyaux $N_k$ sont deux à deux en relation ${R_{\cal P}}$.

Comme $K(\G)$ est par définition le noyau de $\t$, pour tout élément $\p $ de $ K(\G)$ et tout $T\in {\cal C}$, on a $\p([T]) = [T]$. Réciproquement, il s'agit de montrer que, si $\p \in Aut([S])$ satisfait à $\p([T]) = [T]$ pour tout $T\in {\cal C}$, alors $\p([S]_c) = [S]_c$. Or si $N$ est un noyau de $S$ et si $X \in [N]$, on a $\p(X) \in [Cel(N)] \subseteq [S]_c$ ; comme on sait que deux noyaux $N$ et $P$ sont en relation ${R_{\cal P}}$ si et seulement si les cellules $Cel(N)$ et $Cel(P)$ le sont (cf. proposition \ref{Cellules et noyaux d'un relation régulier}), la caractérisation de $[S]_c$ obtenue ci-dessus permet de conclure que l'on a $\p([S]_c) \subseteq [S]_c$. Le même raisonnement appliqué à $\p^{-1}$ nous fournit $\p^{-1}([S]_c)\subseteq [S]_c$ et on a donc le résultat. \qed\\

\begin{nota}[Le sous-groupe {\boldmath $D(\G)$}]\label{def le groupe D}Posons 
$$D(\G) = \{\p \in Aut([S]) \mid \forall ~N\in {\cal N}, \ \p([N]) = [N]\}.$$ 
Comme $[S] = \oplus_{N\in {\cal N}}[N]$, il est clair que $D(\G)$ s'identifie à $\prod_{N\in {\cal N}}Aut([N])$, via l'isomorphisme $\rho_2 : \p \longmapsto (\p_{|[N]})_{N\in {\cal N}}$.\\
\end{nota}

\begin{propo}Le sous-groupe $D(\G)$ de $Aut([S])$ est inclus dans $K(\G)$.
\end{propo}
\proof D'après la proposition précédente, il suffit de vérifier que, pour tout $\p \in D(\G)$ et tout $T \in {\cal C}$, on a $\p([T]) = [T]$. Or toute cellule $T$ est la réunion des $N(s)$, pour $s\in T$ (remarque \ref{props equiv} c), donc $[T] = \sum_{s\in T}[N(s)]$ et on en déduit le résultat.\qed \\

\begin{nota}[Le sous-groupe {\boldmath $K^\circ(\G)$}]\label{Le groupe Krond}Soit $T \in {\cal C}$. On rappelle que $V_T$ désigne l'espace vectoriel quotient $[T]/[T^\star]$. D'après la proposition \ref{sous-espaces [T] permutés}, tout élément $\p$ de $K(\G)$ induit un automorphisme $\bar \p_T$ de l'espace vectoriel quotient $V_T$, donné par $\bar \p_T : \overline X = X + [T^\star] \longmapsto \overline {\p(X)} = \p(X) + [T^\star]$. 

L'application $\rho_3 : \p \longmapsto (\bar\p_T)_{T\in {\cal C}}$ est clairement un morphisme de groupes de $K(\G)$ dans $\prod_{T\in {\cal C}}Aut(V_T)$. On note $K^\circ(\G)$ son noyau. On a $$K^\circ(\G) = \{\p \in K(\G) \mid \forall ~s\in S, \ \p(\{s\}) \in \{s\} + [(C^2(s))^\star]\}.$$
\end{nota}

\begin{propo}\label{decomposition K(Gamma)}La restriction de $\rho_3$ à $D(\G)$ est bijective. On a donc en particulier la décomposition $K(\G) = K^\circ(\G) \rtimes D(\G)$.
\end{propo}
\proof Soient $N\in {\cal N}$ et $T = Cel(N)$. Comme $T$ est la réunion disjointe de $N$ et de $T^\star$, on a la décomposition $[T] = [N] \oplus [T^\star]$, et le morphisme $X \longmapsto \overline X = X + [T^\star]$ est donc un isomorphisme de $[N]$ sur $V_T$. Si l'on identifie $[N]$ à $V_T$ (et $Aut([N])$ à $Aut(V_T)$) via cet isomorphisme, alors la restriction de $\rho_3$ à $D(\G)$ s'identifie à l'isomorphisme $\rho_2 : \p \longmapsto (\p_{|[N]})_{N\in {\cal N}}$ de $D(\G)$ sur $\prod_{N\in {\cal N}}Aut([N])$ (cf. notation \ref{def le groupe D}). On a donc le résultat. \qed \\

\begin{comms}\label{exemples d'éléments de D et Kzéro}{\bf a.} On a $K(\G) \cap Aut(\G) \subseteq D(\G)$. De plus, l'isomorphisme $\rho_1$ de la proposition \ref{intersection de Aut(G) et K} est induit par l'isomorphisme $\rho_2$ si l'on identifie, pour tout $N\in {\cal N}$, $Perm(N)$ à un sous-groupe de $Aut([N])$, comme en remarque \ref{Rq tilde f} b.

{\bf b.} Soient $s \in S$ et $t \in C^2(s)$, $t\not= s$. Alors l'élément $\a_{s,t}$ de $Aut([S],[S]_c)$, défini en notation \ref{automs élémentaires} (cf. remarque \ref{rems automs éléms} b), appartient à $D(\G)$ ou à $K^\circ(\G)$ selon que $C^2(t) = C^2(s)$ ou que $C^2(t) \varsubsetneq C^2(s)$.
\end{comms}

\subsection{\'Etude du groupe {\boldmath $K^\circ(\G)$}.}\label{K0}

\subsubsection{Sous-groupes de {\boldmath $K^\circ(\G)$}.}\label{sous-groupes de Kzéro}

Les résultats de cette section \ref{sous-groupes de Kzéro} sont valables pour le groupe $K(\G)$ (en remplaçant systématiquement $K^\circ$ par $K$ dans les énoncés qui suivent). Cependant, comme ils ne nous serviront qu'à décrire le groupe $K^\circ(\G)$, c'est dans ce cadre que nous les présentons. \\
 
\begin{defi}[Support]Soit $End([S])$ l'ensemble des endomorphismes de l'espace vectoriel $[S]$. Pour $\p \in End([S])$, on appelle {\it support de $\p$} l'ensemble $D_\p = \{s \in S\mid \p(\{s\}) \not= \{s\}\}$. \\
\end{defi}

\begin{rema}\label{rems support}{\bf a.} Soit $\p\in End([S])$. On vérifie facilement que l'on a $D_\p = \emptyset \Longleftrightarrow \p = Id_{[S]}$, et que, si $\p = \p_1\circ \cdots \circ \p_n$, avec $\p_1, \ldots ,~\p_n \in End([S])$, alors $D_{\p} \subseteq \bigcup_{1\leq k \leq n}D_{\p_k}$. 

{\bf b.} De plus, si $\p \in Aut([S])$, alors $D_{\p^{-1}} = D_\p$.\\
\end{rema}

\begin{nota}[Les sous-groupes {\boldmath $K^\circ_Y(\G)$}]Soit $Y \subseteq S$. Nous notons $K^\circ_Y(\G)$ l'ensemble $\{\p\in K^\circ(\G) \mid D_\p \subseteq Y\}$. D'après les remarques \ref{rems support} a et b, $K^\circ_Y(\G)$ est un sous-groupe de $K^\circ(\G)$.

Remarquons que $K^\circ_S(\G) = K^\circ(\G)$.\\
\end{nota}

Dans la suite de cette section, nous voulons décrire le groupe $K^\circ_Y(\G)$, pour $Y\subseteq S$, à partir de certains de ses sous-groupes $K^\circ_Z(\G)$, où $Z\subseteq Y$. \\

\begin{defi}[Parties saturées]Soit $Y \subseteq S$. Nous disons qu'une partie $X$ de $S$ est {\it $Y$-saturée} si, pour tout $x \in X\cap Y$, $C^2(x) \cap Y \subseteq X$. Nous disons simplement {\it saturée} pour $S$-saturée.\\
\end{defi}

\begin{rema}\label{propriétés parties saturées}Soient $X$ et $Y$ deux parties de $S$.

{\bf a.} $X$ est $Y$-saturée si et seulement si $X\cap Y = \bigcup_{x \in X\cap Y}C^2(x)\cap Y$. En particulier, les parties saturées sont les réunions de cellules.

{\bf b.} Si $Y'\subseteq Y$ et si $X$ est $Y$-saturée, alors $X$ est $Y'$-saturée.\\
\end{rema}

\begin{defi}[Troncature]Soient $X \subseteq S$ et $\p \in End([S])$. On appelle {\it troncature de $\p$ suivant $X$} l'endomorphisme $\p_X$ de $[S]$ donné par $\p_X(\{s\}) = \p(\{s\})$ si $s\in X$ et $\p_X(\{s\}) = \{s\}$ si $s\not \in X$. Clairement, $\p_X$ coïncide avec $\p$ sur $[X \cup (S\setminus D_\p)]$ et $D_{\p_X} = X\cap D_\p$.\\
\end{defi}

\begin{lem}\label{troncature}Soient $Y \subseteq S$ et $X$ une partie $Y$-saturée. 
\begin{enumerate}
\item Soient $\p,~\psi \in K^\circ_Y(\G)$. Alors $(\p\circ \psi)_X = \p_X \circ \psi_X$.
\item L'application $K^\circ_Y(\G) \longrightarrow K^\circ_{X\cap Y}(\G)$, $\p \longmapsto \p_X$ est un morphisme de groupes de noyau $K^\circ_{Y\setminus X}(\G)$. 
\end{enumerate}
\end{lem}
\proof Montrons le premier point. Soient $\p,~\psi \in K^\circ_Y(\G)$ et $s \in S$. Montrons que $(\p \circ \psi)_X(\{s\}) = (\p_X \circ \psi_X)(\{s\})$. Si $s \not \in X\cap Y$, alors $\p_X(\{s\}) = \psi_X(\{s\}) = (\p \circ \psi)_X(\{s\}) = \{s\}$ et le résultat est clair. Si $s \in X\cap Y$, alors $C^2(s) \subseteq X \cup (S\setminus Y)$ puisque $X$ est $Y$-saturée. On a alors $C^2(s) \subseteq X \cup (S\setminus D_\p)$ (puisque $D_\p \subseteq Y$) et $\p_X$ coïncide donc avec $\p$ sur $[C^2(s)]$. De plus, comme $\psi \in K(\G)$, on a $\psi_X(\{s\}) = \psi(\{s\}) \in [(C^2(s)]$, d'où l'on déduit que $(\p_X \circ \psi_X) (\{s\}) = \p(\psi(\{s\})) = (\p \circ \psi)_X(\{s\})$.

Montrons le second point. Pour montrer que $\p \longmapsto \p_X$ est un morphisme de groupes de $K^\circ_Y(\G)$ dans $K^\circ_{X\cap Y}(\G)$, il suffit (grâce au premier point) de montrer que, pour $\p \in K^\circ_Y(\G)$, la troncature $\p_X$ appartient à $K^\circ_{X\cap Y}(\G)$. Vu la définition de $\p_X$, le seul fait non trivial à montrer est que $\p_X$ appartient à $K(\G)$. Le premier point nous montre en particulier que $(\p^{-1})_X$ et $\p_X$ sont des automorphismes (de $[S]$) inverses l'un de l'autre (puisque $(Id_{[S]})_X = Id_{[S]}$). D'après la proposition \ref{caract de F et K}, il suffit donc de montrer que l'on a $\p_X([T]) = [T]$ pour toute cellule $T$. Or, par définition de $\p_X$ et $(\p^{-1})_X$ et comme $\p \in K(\G)$, on a clairement $\p_X([T]) \subseteq [T]$ et $(\p^{-1})_X([T]) = (\p_X)^{-1}([T]) \subseteq [T]$, d'où le résultat. Le noyau de ce morphisme est clairement le sous-groupe $K^\circ_{Y\setminus X}(\G)$ de $K^\circ_Y(\G)$, et le lemme est donc démontré.\qed \\

\begin{propo}\label{decomp}Soit $Y \subseteq S$ et soit $X$ une partie $Y$-saturée de $S$. On a la décomposition : $$K^\circ_Y(\G) = K^\circ_{Y\setminus X}(\G) \rtimes K^\circ_{X \cap Y}(\G).$$
\end{propo}
\proof L'assertion {\it 2} du lemme précédent nous dit que l'on a un morphisme de groupes $\p \longmapsto \p_X$ de $K^\circ_Y(\G)$ dans $K^\circ_{X\cap Y}(\G)$ et que ce morphisme a pour noyau $K^\circ_{Y\setminus X}(\G)$. Comme l'inclusion naturelle $K^\circ_{X\cap Y}(\G) \hookrightarrow K^\circ_Y(\G)$ en est clairement une section, on a le résultat. \qed \\

Ce résultat suggère une méthode pour décrire le sous-groupe $K^\circ_Y(\G)$, pour une partie donnée $Y$ de $S$ : trouver une partie $Y$-saturée $X \varsubsetneq Y$ telle que le sous-groupe $K^\circ_{Y\setminus X}(\G)$ soit "simple à décrire", puis décrire par récurrence le sous-groupe $K^\circ_{X\cap Y}(\G) = K^\circ_{X}(\G)$.

Dans la section suivante, nous décrivons le sous-groupe $K^\circ_Z(\G)$ de $K^\circ(\G)$ dans le cas où $Z$ est une partie d'"épaisseur nulle" (voir la définition \ref{def parties de profondeur nulle} ci-dessous).

Dans la section \ref{profondeur}, nous allons voir comment décomposer toute partie $Y$ de $S$ en l'union disjointe d'une partie d'épaisseur nulle $Y_0$ et d'une partie $Y$-saturée $Y_{\geq 1}$. La proposition \ref{decomp}, appliquée récursivement, nous permettra alors de décomposer le groupe $K^\circ_Y(\G)$ (au moins lorsque $Y$ est "d'épaisseur finie", voir la définition \ref{def profondeur finie} ci-dessous) en produits semi-directs itérés de certains de ses sous-groupes $K^\circ_Z(\G)$ avec $Z$ d'épaisseur nulle (théorème \ref{décomp de KY} ci-dessous).

\subsubsection{Parties d'épaisseur nulle.}\label{Parties de profondeur nulle.}

Nous définissons plus loin (en définition \ref{def profondeur finie}) l'"épaisseur" (dans $\N \cup \{\infty\}$) d'une partie de $S$. Définissons pour le moment les parties d'"épaisseur nulle" : \\ 

\begin{defi}[Parties d'épaisseur nulle]\label{def parties de profondeur nulle}On dit qu'une partie $Z$ de $S$ est {\it d'épaisseur nulle} si, pour tout $s\in Z$, $(C^2(s))^\star \cap Z = \emptyset$ (autrement dit si tous les noyaux sont $Z$-saturés).\\
\end{defi}

\begin{rema}\label{rems épaisseur nulle}{\bf a.} Si $Z$ est d'épaisseur nulle et si $Z'\subseteq Z$, alors $Z'$ est d'épaisseur nulle.

{\bf b.} Les noyaux (et donc les parties de noyaux, d'après a) sont d'épaisseur nulle.\\
\end{rema}

Fixons une partie $Z$ de $S$ d'épaisseur nulle.\\

\begin{lem}\label{lemme technique}Soient $\p, ~\psi \in End([S])$ de supports inclus dans $Z$ et soit $s \in Z$. Si on a $\p(\{s\}) = \{s\} + Y$ et $\psi(\{s\}) = \{s\} + Y'$, où $Y,~Y' \in [(C^2(s))^\star]$, alors $(\p\circ \psi)(\{s\}) = (\psi \circ \p)(\{s\}) = \{s\} + Y + Y'$.
\end{lem}
\proof Par hypothèse, $Z$ est d'épaisseur nulle et contient $s$, donc est disjointe de $(C^2(s))^\star$. A fortiori, les parties $Y$ et $Y'$ de $(C^2(s))^\star$ sont disjointes de $D_\p$ et de $D_\psi$ (qui sont inclus dans $Z$). On a donc $\p(Y') = Y'$, $\psi(Y) = Y$, et un calcul direct nous donne le résultat.\qed \\

\begin{propo}\label{caracK0de profondeur nulle}$K^\circ_Z(\G) = \{ \p \in End([S]) \mid D_\p \subseteq  Z \textrm{ et, } \forall ~s \in Z, \ \p(\{s\}) \in \{s\} + [(C^2(s))^\star]\}$. En particulier, $K^\circ_Z(\G)$ est un 2-groupe élémentaire.
\end{propo}
\proof Par définition de $K^\circ(\G)$ et de son sous-groupe $K^\circ_Z(\G)$, l'inclusion $K^\circ_Z(\G) \subseteq  \{ \p \in End([S]) \mid D_\p \subseteq  Z \textrm{ et, } \forall ~s \in Z, \ \p(\{s\}) \in \{s\} + [(C^2(s))^\star]\}$ est évidente. Soit donc $\p$ un endomorphisme de $[S]$ tel que $D_\p \subseteq  Z$ et que $\p(\{s\}) \in \{s\} + [(C^2(s))^\star]$, pour tout $s \in Z$.

Montrons que $\p^2 = Id_{[S]}$. Si $s\not \in Z$, alors on a $\p(\{s\}) = \{s\}$, donc $\p^2(\{s\}) = \{s\}$ ; si $s\in Z$, alors $\p(\{s\}) = \{s\} + Y$, où $Y \in [(C^2(s))^\star]$, et le lemme précédent nous fournit $\p^2(\{s\}) = \{s\} + Y + Y = \{s\}$.

Il reste à montrer que $\p$ appartient à $K(\G)$. Il suffit pour cela de montrer que, pour toute cellule $T$, $\p([T]) = [T]$ (cf. proposition \ref{caract de F et K}) et, comme $\p$ est involutif, il suffit de montrer que $\p([T]) \subseteq [T]$. Or pour $s \in T$ (qui satisfait donc à $C^2(s)\subseteq T$), on a soit $s \not \in D_\p$, auquel cas $\p(\{s\}) = \{s\} \in [T]$, soit $s \in D_\p$, donc $s \in Z$ et $\p(\{s\}) \in \{s\} + [(C^2(s))^\star] \subseteq [T]$. On a donc le résultat.\qed \\

\begin{cor}\label{proddirectétal}Soit $(Z_i)_{i \in I}$ une partition de $Z$. Alors l'application $\p \longmapsto (\p_{Z_i})_{i\in I}$, où $\p_{Z_i}$ est la troncature de $\p$ suivant $Z_i$, est un isomorphisme de groupes de $K^\circ_Z(\G)$ sur le produit direct de ses sous-groupes $K^\circ_{Z_i}(\G) \ (i\in I)$.
\end{cor}
\proof Remarquons que, puisque $Z$ est d'épaisseur nulle, chaque $Z_i$, $i\in I$, est d'épaisseur nulle (cf. remarque \ref{rems épaisseur nulle} a). D'après la caractérisation des éléments de $K^\circ_X(\G)$ pour une partie d'épaisseur nulle $X$, obtenue dans la proposition \ref{caracK0de profondeur nulle}, il est clair que, pour tout $i \in I$, la troncature de tout élément de $K^\circ_Z(\G)$ suivant $Z_i$ est un élément de $K^\circ_{Z_i}(\G)$. 

De plus, si $\p, ~\psi \in K^\circ_Z(\G)$, alors $(\p \circ \psi)_{Z_i} = \p_{Z_i}\circ \psi_{Z_i}$. En effet, pour $s \in  Z_i$, on a $\p_{Z_i}(\{s\}) = \p(\{s\}) = \{s\} + Y$ et $\psi_{Z_i}(\{s\}) = \psi(\{s\}) = \{s\} + Y'$, où $Y,~ Y' \in [(C^2(s))^\star]$, et le lemme \ref{lemme technique} nous fournit $(\p \circ \psi)(\{s\}) = \{s\} + Y + Y' = (\p_{Z_i} \circ \psi_{Z_i})(\{s\})$. L'application $\p \longmapsto (\p_{Z_i})_{i\in I}$ est donc un morphisme de groupes de $K^\circ_Z(\G)$ dans $\prod_{i\in I}K^\circ_{Z_i}(\G)$. Ce morphisme est injectif car l'image $(\p_{Z_i})_{i\in I}$ détermine $\p$ sur $[Z]$ (puisque $(Z_i)_{ i\in I}$ est une partition de $Z$), donc sur $[D_\p]$. 

Montrons qu'il est surjectif. Si $(\p_{i})_{i\in I}$ est un élément de $\prod_{i\in I}K^\circ_{Z_i}(\G)$, définissons l'endomorphisme $\p$ de $[S]$ par $\p(\{s\}) = \{s\}$ si $s\not \in Z$, et $\p(\{s\}) = \p_i(\{s\})$ si $s\in Z_i$. D'après la caractérisation de la proposition \ref{caracK0de profondeur nulle}, $\p$ est un élément de $K^\circ_Z(\G)$, et il s'envoie clairement sur la famille $(\p_{i})_{i\in I}$. On a donc le résultat.  \qed \\

\begin{cor}\label{K0Xpour Xnoyau}Soient $N \in {\cal N}$, et $X \subseteq N$. On pose $T = Cel(N)$. Alors $X$ est d'épaisseur nulle et l'application $\p \longmapsto (Id_{[S]} + \p)_{|[X]}$ est un isomorphisme de 2-groupes élémentaires de $K^\circ_X(\G)$ sur le groupe (additif) ${\cal L}_{\F_2}([X],[T^\star])$, constitué des applications $\F_2$-linéaires de $[X]$ dans $[T^\star]$.
\end{cor}
\proof La partie $X$ de $N$ est d'épaisseur nulle d'après la remarque \ref{rems épaisseur nulle} b. Soit $\p \in K^\circ_X(\G)$. On sait que, pour $s\in S$, on a $\p(\{s\}) \in \{s\} + [(C^2(s))^\star]$. En particulier, pour $s\in X \subseteq N$, on a $C^2(s) = T$, donc $\p(\{s\}) \in \{s\} + [T^\star]$. L'application $\F_2$-linéaire $(Id_{[S]} + \p)_{|[X]}$ est donc à valeurs dans $[T^\star]$ et l'application $\p \longmapsto (Id_{[S]} + \p)_{|[X]}$ est donc bien à valeurs dans ${\cal L}_{\F_2}([X],[T^\star])$. 

Soient $\p, ~\psi \in K^\circ_X(\G)$ et $s \in X$. On a $\p(\{s\}) = \{s\} + Y$ et $\psi(\{s\}) = \{s\} + Y'$, où $Y,~ Y' \in [T^\star]$, et, d'après le lemme \ref{lemme technique}, $\p \circ \psi(\{s\}) = \{s\} + Y + Y'$, d'où $(Id_{[S]} + \p \circ \psi)(\{s\}) =  Y + Y' = (Id_{[S]} + \p)(\{s\}) + (Id_{[S]} + \psi)(\{s\})$. Donc $\p \longmapsto (Id_{[S]} + \p)_{|[X]}$ est un morphisme de groupes. Ce morphisme est injectif, car l'image $(Id_{[S]} + \p)_{|[X]}$ de $\p$ détermine $\p$ sur $[X]$, donc sur $[D_\p]$. 

Montrons qu'il est surjectif. Pour $f \in {\cal L}_{\F_2}([X],[T^\star])$, notons encore $f$ l'endomorphisme $f \oplus 0_{|[S\setminus X]}$ (prolongement de $f$ à $[S]$, défini par $0$ sur $[S\setminus X]$). Alors l'endomorphisme $\p = Id_{[S]} + f$ de $[S]$ appartient à $K^\circ_X(\G)$ (d'après la proposition \ref{caracK0de profondeur nulle}), et s'envoie clairement sur $f$. \qed \\

Pour décomposer $K^\circ_Z(\G)$ (où $Z$ est d'épaisseur nulle) comme dans le corollaire \ref{proddirectétal}, nous allons privilégier, vu le corollaire \ref{K0Xpour Xnoyau}, les partitions de $Z$ en parties de noyaux. Par exemple, on obtient :

\begin{description}
\item[(E1)] $K^\circ_Z(\G) \approx \prod_{N\in {\cal N}(Z)}{\cal L}_{\F_2}([N\cap Z],[(Cel(N))^\star])$, via $\p \mapsto ((Id_{[S]} + \p)_{|[N\cap Z]})_{N\in {\cal N}(Z)}$, où ${\cal N}(Z) = \{N(z) \mid z\in Z\}$ est l'ensemble des noyaux rencontrés par $Z$, \\
\item[(E2)] $K^\circ_Z(\G) \approx \prod_{s\in Z}{\cal L}_{\F_2}([\{s\}],[(C^2(s))^\star]) \approx \prod_{s\in Z}[(C^2(s))^\star]$, via $\p \mapsto ((Id_{[S]} + \p)_{|[\{s\}]})_{s\in Z}$ et $(f_s)_{s\in Z} \mapsto (f_s(\{s\}))_{s\in Z}$, avec la partition de $Z$ en l'ensemble de ses singletons.
\end{description}

\subsubsection{Profondeur et décomposition en parties d'épaisseur nulle.}\label{profondeur}

\begin{defi}[Profondeur]Soit $Y \subseteq S$. On appelle {\it chaîne (d'éléments de $Y$)} toute suite finie $(y_0, ~y_1, \ldots , y_p)$ d'éléments de $Y$ telle que, pour $1\leq i \leq p$, $C^2(y_{i-1}) \varsubsetneq C^2(y_i)$. Si $(y_0, ~y_1, \ldots , y_p)$ est une chaîne, on dit que sa {\it longueur} est $p$, et qu'elle a pour {\it origine} $y_0$.

Soit $Y \subseteq S$. On définit la fonction {\it $Y$-profondeur} $p_Y : Y \longrightarrow\N \cup\{\infty\}$, de la façon suivante : pour $y \in Y$, $p_Y(y)$ est la borne supérieure de l'ensemble des longueurs des chaînes d'éléments de $Y$ d'origine $y$. Pour $p \in \N \cup \{\infty\}$, on pose $Y_p = \{y\in Y \mid p_Y(y) = p\}$, et $Y_{\geq p} = \cup_{k\geq p}Y_k = \{y\in Y \mid p_Y(y) \geq p \}$ (En particulier, $Y_{\geq 0} = Y$ et $Y_\infty = Y_{\geq \infty}$).\\
\end{defi}

\begin{rema}\label{propriété profondeur}{\bf a.} Si $Y$ est une cellule, alors $Y_0 = Noy(Y)$, et $Y_{\geq 1} = Y \setminus Noy(Y) = Y^\star$. Plus généralement, si $Y$ est saturée (i.e. si $Y$ est une réunion de cellules) et si $p \in \N \cup \{\infty\}$, alors $Y_p$ est une réunion de noyaux, et $Y_{\geq p}$ est saturée.

{\bf b.} Les $Y_p$, pour $p \in \N \cup \{\infty\}$, sont deux à deux disjoints, et chaque partie $Y_{\geq p}$ est $Y$-saturée.

{\bf c.} Soit $n \in \N$. Alors $Y_n$ est d'épaisseur nulle, et satisfait à $Y_n \not = \emptyset \Longrightarrow Y_k \not = \emptyset$, pour tout $k \leq n$. De plus, si $p \in \N\cup \{\infty\}$, on vérifie que $(Y_{\geq n})_p = Y_{n+p}$, d'où $(Y_{\geq n})_{\geq p} = Y_{\geq n+p}$.\\
\end{rema}

\begin{defi}[\'Epaisseur]\label{def profondeur finie}Soit $Y\subseteq S$. Nous appelons {\it épaisseur de $Y$} la borne supérieure (dans $\N\cup\{\infty\}$) de la fonction $p_Y$. 

On vérifie facilement que les parties d'épaisseur 0, au sens que l'on vient de définir, sont les parties d'épaisseur nulle au sens de la définition \ref{def parties de profondeur nulle} de la section \ref{Parties de profondeur nulle.}.


Lorsque $S$ est d'épaisseur finie, nous disons que $(W,S)$ et $\G$ sont {\it d'épaisseur finie}. Comme cette notion ne dépend pas de $S\in {\cal S}(W)$ (puisque $W$ est rigide), nous disons aussi, sans ambiguïté, que $W$ est d'épaisseur finie. C'est par exemple le cas lorsque ${\cal C}$ (ou ${\cal N}$) est fini et, en particulier, on voit donc que tout groupe de Coxeter de rang fini est d'épaisseur finie.\\
\end{defi}

\begin{theo}\label{décomp de KY}Soit $Y \subseteq S$ d'épaisseur finie $e$. Alors :
$$K^\circ_Y(\G) = K^\circ_{Y_0}(\G) \rtimes (K^\circ_{Y_1}(\G) \rtimes (\cdots \rtimes (K^\circ_{Y_{e-1}}(\G) \rtimes K^\circ_{Y_{e}}(\G))\cdots)).$$

De plus, si $Y$ est saturée, alors $K^\circ_{Y_{e}}(\G) = \{Id_{[S]}\}$.
\end{theo}
\proof Soit $k \in \N$. Sans hypothèse sur $Y$, la proposition \ref{decomp}, appliquée à $Y_{\geq k}$ et à la partie $Y_{\geq k}$-saturée $(Y_{\geq k})_{\geq 1} = Y_{\geq k+1}$ (remarque \ref{propriété profondeur} c), nous fournit la décomposition : 
\begin{equation*}K^\circ_{Y_{\geq k}}(\G) = K^\circ_{Y_k}(\G) \rtimes K^\circ_{Y_{\geq k+1}}(\G).
\end{equation*}
On suppose ici $Y$ d'épaisseur finie $e$, donc on a $Y_{\geq e+1} = \emptyset$ et $K^\circ_{Y_{\geq e+1}}(\G) = \{Id_{[S]}\}$. La formule voulue s'obtient alors facilement par récurrence, en appliquant le procédé ci-dessus successivement aux parties $Y_{\geq k}$, pour $0\leq k \leq e$. 

Si, de plus, $Y$ est saturée, alors, pour tout $n \in \N$, on a $y \in Y_n \Rightarrow (C^2(y))^\star \subseteq Y_{\geq n+1}$, et on voit que $Y_e$ est nécessairement constitué d'éléments $y$ de $Y$ tels que $(C^2(y))^\star = \emptyset$. Mais alors la formule (E2) appliquée à $Z = Y_e$ (qui est d'épaisseur nulle, d'après la remarque \ref{propriété profondeur} c) nous montre que l'on a $K^\circ_{Y_{e}}(\G) = \{Id_{[S]}\}$. \qed 

\subsection{Conclusion.}\label{conclusion}

On rappelle que nous voulions, dans cette partie, décrire le groupe $Aut([S],[S]_c)$.\\

Nous avons défini (en section \ref{premier dévissage}) un sous-groupe $G$ de $Aut(\G)$ isomorphe à $Aut({\cal N},{R_{\cal P}},Card)$ pour lequel on a (d'après la proposition \ref{premières décompositions}) : 
\begin{equation}\label{1}
Aut([S],[S]_c) = K(\G) \rtimes G \textrm{   et  } Aut(\G) = (K(\G)\cap Aut(\G)) \rtimes G.
\end{equation}

Les propositions \ref{intersection de Aut(G) et K} (section \ref{premier dévissage}) et \ref{decomposition K(Gamma)} (section \ref{le groupe K}) nous donnent (avec le commentaire \ref{exemples d'éléments de D et Kzéro} a): 
\begin{equation}\label{2}
\begin{array}{c}
K(\G) = K^\circ(\G) \rtimes D(\G), \textrm{  avec } D(\G) \approx  \prod_{N \in {\cal N}}Aut([N]), \textrm{ via } \rho_2 : \p \mapsto (\p_{|[N]})_{N\in {\cal N}},\\
K(\G)\cap Aut(\G)\subseteq D(\G), \textrm{ et } K(\G)\cap Aut(\G) \approx  \prod_{N \in {\cal N}}Perm(N), \textrm{ via } \rho_1 : \s \mapsto (\s_{|N})_{N\in {\cal N}}.
\end{array}
\end{equation}

Le théorème \ref{décomp de KY} (section \ref{profondeur}) nous permet de décrire le groupe $K^\circ(\G) = K^\circ_S(\G)$, dans le cas où $S$ (qui est une partie saturée) est d'épaisseur finie $e$. On obtient :
\begin{equation}\label{3}
K^\circ(\G) = K^\circ_{S_0}(\G) \rtimes (K^\circ_{S_1}(\G) \rtimes (\cdots \rtimes (K^\circ_{S_{e-2}}(\G) \rtimes K^\circ_{S_{e-1}}(\G))\cdots)).
\end{equation}

Enfin, pour $0\leq k \leq e-1$, la partie $S_k$ de $S$ est d'épaisseur nulle (cf. remarque \ref{propriété profondeur} c) et contient les noyaux qu'elle rencontre (cf. remarque \ref{propriété profondeur} a), ce qui signifie que, ${\cal N}(S_k)$ désignant l'ensemble $\{N(s) \mid s \in S_k\}$, on a, pour tout $N \in {\cal N}(S_k)$, $N \subseteq S_k$ ; les formules (E1) et (E2) appliquées à $S_k$ nous donnent donc : 
\begin{equation}\label{4}
K^\circ_{S_k}(\G) \approx \begin{cases}
\prod_{N\in {\cal N}(S_k)}{\cal L}_{\F_2}([N],[(Cel(N))^\star]), \textrm{  via } \p \mapsto ((Id_{[S]} + \p)_{|[N]})_{N\in {\cal N}(S_k)}\\
\prod_{s\in S_k}[(C^2(s))^\star], \textrm{  via } \p \mapsto (\{s\} + \p(\{s\}))_{s\in S_k}.\end{cases}
\end{equation}\\

\begin{rema}\label{remarque finale}Rappelons que l'on note ${\cal E}(S)$ l'ensemble des singletons de $S$. On a : 

{\bf a.} $S$ est d'épaisseur nulle $\Leftrightarrow {\cal N} = {\cal C} \Leftrightarrow K^\circ(\G) = \{1\}$,

{\bf b.} $ {\cal N} = {\cal E}(S) \Leftrightarrow D(\G) = \{1\} \Leftrightarrow Aut(\G)\cap K(\G) = \{1\} \Leftrightarrow G = Aut(\G)$,

{\bf c.} $\G$ satisfait aux conditions de a et de b $\Leftrightarrow {\cal C} = {\cal E}(S) \Leftrightarrow K(\G) = \{1\} \Leftrightarrow Aut([S],[S]_c) = Aut(\G)$. On retrouve ainsi le résultat de la proposition \ref{caract}.
\end{rema}

\section*{Exemples.}

Dans les exemples de graphe de Coxeter qui suivent, toutes les arêtes sont étiquetées $\infty$.\\


\begin{ex}Nous retrouvons les résultats obtenus par J. Tits dans la partie finale de [T] (le groupe $Aut(F(\G))$ étudié par J. Tits est isomorphe à $Aut([S],[S]_c)$).
\begin{itemize}
\item[$\bullet$] Supposons que $S$ soit la réunion disjointe de deux sous-ensembles non vides $S'$ et $S''$ tels que $m_{s,t} = \infty \Leftrightarrow (s,t) \in S' \times S''$ ou $(s,t) \in S''\times S'$ (on dit alors que le graphe de $\G$ est "bipartite complet"). Avec notre terminologie, il s'agit du cas où $S$ est d'épaisseur nulle et a deux noyaux (qui sont $S'$ et $S''$). 

On a alors, d'après les formules (\ref{1}) et (\ref{2}) et la remarque \ref{remarque finale} a, $Aut([S],[S]_c) = D(\G) \rtimes G \approx (Aut([S'])\times Aut([S'']))\rtimes G$. De plus, on a clairement $G \approx Aut({\cal N},{R_{\cal P}},Card) = \{1\}$ si les noyaux $S'$ et $S''$ n'ont pas même cardinal, et $G \approx \Z/2\Z$ si $S'$ et $S''$ ont même cardinal.
\item[$\bullet$] Supposons que le graphe de $\G$ soit un cycle de longueur $n \geq 5$. On vérifie facilement que ${\cal C} = {\cal E}(S)$ et donc $Aut([S],[S]_c) = Aut(\G)$ (cf. proposition \ref{caract}, ou remarque \ref{remarque finale} c).
\item[$\bullet$] Les autres cas étudiés par J. Tits sont d'épaisseur 1, et satisfont à ${\cal N} = {\cal E}(S)$. On a alors $Aut([S],[S]_c) = K^\circ_{S_0}(\G) \rtimes Aut(\G)$ (cf. formules (\ref{1}), (\ref{2}), (\ref{3}) et remarque \ref{remarque finale} b) et  $K^\circ_{S_0}(\G) \approx \prod_{s\in S_0}[(C^2(s))^\star]$ (cf. formule (\ref{4})).\\
\end{itemize}
\end{ex}

\begin{ex}Soit $\G$ donnée par le graphe \begin{picture}(75,20)(-25,-2)
\put(0,0){\circle{4}}\put(2,4){3}
\put(20,0){\circle{4}}\put(22,4){4}
\put(40,0){\circle{4}}\put(42,4){5}
\put(2,0){\line(1,0){16}}
\put(22,0){\line(1,0){16}}
\put(2,0){\line(1,0){16}}
\put(-2,-1){\line(-3,-2){12}}
\put(-15,-10){\circle{4}}\put(-25,-10){2}
\put(-2,1){\line(-3,2){12}}
\put(-15,10){\circle{4}}\put(-25,10){1}
\end{picture}. On a :

$$\begin{array}{|l|c|c|c|c|c|}\hline
\textrm{Sommet $s$} & 1 & 2 & 3 & 4 & 5 \\ \hline
\textrm{Cellule $C^2(s)$} & \{1,2\} & \{1,2\} & \{3,5\} & \{1,2,4\} & \{5\} \\ \hline
\textrm{Noyau $N(s)$} & \{1,2\} & \{1,2\} & \{3\} & \{4\} & \{5\}  \\ \hline 
\end{array}$$

On vérifie que $G \approx Aut({\cal N},{R_{\cal P}},Card) = \{1\}$ et que $S$ est d'épaisseur 1, avec $S_0 = \{3,4\}$ et $S_1 = \{1,2,5\}$. On a alors, grâce aux formules (1) à (4),  $Aut([S],[S]_c) = K^\circ_{S_0}(\G) \rtimes D(\G)$, avec :  
$$D(\G) \approx Aut([\{1,2\}]) \approx GL_2(\F_2) \ \textrm{ et }  K^\circ_{S_0}(\G)\approx [(C^2(3))^\star]\times [(C^2(4))^\star] \approx (\F_2)^3.$$
\end{ex}

\vspace{0.5cm}

\begin{ex}Soit $\G$ donnée par le graphe infini (indexé par $\Z$) : 
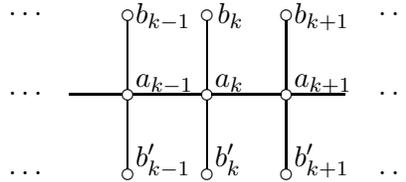
\begin{center}
\begin{picture}(80,80)(-10,-10)
\put(-45,30){$\ldots$}
\put(-45,0){$\ldots$}
\put(-45,60){$\ldots$}
\put(95,30){$\ldots$}
\put(95,0){$\ldots$}
\put(95,60){$\ldots$}

\put(-2,30){\line(-1,0){20}}
\put(62,30){\line(1,0){20}}

\put(0,0){\circle{4}}\put(3,3){$b'_{k-1}$}
\put(30,0){\circle{4}}\put(33,3){$b'_k$}
\put(60,0){\circle{4}}\put(63,3){$b'_{k+1}$}
\put(0,30){\circle{4}}\put(3,33){$a_{k-1}$}
\put(30,30){\circle{4}}\put(33,33){$a_k$}
\put(60,30){\circle{4}}\put(63,33){$a_{k+1}$}
\put(0,60){\circle{4}}\put(3,57){$b_{k-1}$}
\put(30,60){\circle{4}}\put(34,57){$b_k$}
\put(60,60){\circle{4}}\put(63,57){$b_{k+1}$}

\put(0,2){\line(0,1){26}}
\put(0,32){\line(0,1){26}}
\put(30,2){\line(0,1){26}}
\put(30,32){\line(0,1){26}}
\put(60,2){\line(0,1){26}}
\put(60,32){\line(0,1){26}}

\put(2,30){\line(1,0){26}}
\put(32,30){\line(1,0){26}}
\end{picture}
\end{center}

On a : 
$$\begin{array}{|l|c|c|c|}\hline
\textrm{Sommet $s$} & a_k & b_k & b_k' \\ \hline
\textrm{Cellule $C^2(s)$} & \{a_k, b_{k-1},b'_{k-1}, b_{k+1},b'_{k+1}\} & \{b_k,b'_k\} & \{b_k,b'_k\}  \\ \hline
\textrm{Noyau $N(s)$} & \{a_k\} & \{b_k,b'_k\} & \{b_k,b'_k\}  \\ \hline 
\end{array}$$

Il est clair que tout élément de $Aut({\cal N},{R_{\cal P}},Card)$ est entièrement déterminé par son action sur l'ensemble de noyaux $\{\{a_k\}\mid k\in \Z\} \approx \Z$, puisque si $\{a_k\}$ s'envoie sur $\{a_l\}$, alors nécessairement, $\{b_k,b_k'\}$ s'envoie sur $\{b_l,b'_l\}$. On voit alors que $Aut({\cal N},{R_{\cal P}},Card)$ s'identifie au groupe diédral infini $D_\infty = <\psi,~\psi'>$ engendré par les symétries $\psi : k\mapsto -k, ~k\in \Z$ et $\psi' : k\mapsto -k+1, ~k\in \Z$. Via cette identification (et celle de $Aut({\cal N},{R_{\cal P}},Card)$ au sous-groupe $G$ de $Aut(\G)$), on a donc $Aut([S],[S]_c) = K(\G) \rtimes D_\infty$ (cf. formule (1)).

On déduit du tableau ci-dessus que $S$ est d'épaisseur 1, avec $S_0 = \{a_k \mid k\in \Z \}$ et $S_1 = \{b_k, b'_k \mid k\in \Z \}$. On obtient, grâce aux formules (\ref{2}) (\ref{3}) et (\ref{4}), $K(\G) = K^\circ_{S_0}(\G) \rtimes D(\G)$, avec :
$$D(\G)\approx \prod_{k\in \Z}GL_2(\F_2)\  \textrm{ et } K^0_{S_0}(\G) \approx \prod_{k\in \Z}(\F_2)^4 \approx (\F_2)^\Z.$$
\end{ex}

\vspace{0.5cm}

\begin{ex}Soit $e \in \N^\star \cup\{\infty\}$ et soit $\G = (m_{s,t})_{s,t\in S}$, où $S = \{s_k,~s_k' \mid 0\leq k < e+1\}$ et $m_{s,t} = \infty \Longleftrightarrow \{s,t\} = \{s_k,s'_l\}$ avec $0 \leq l\leq k <e+1$. On a :
$$\begin{array}{|l|c|c|}\hline
\textrm{Sommet $s$} & s_k & s'_k \\ \hline
\textrm{Cellule $C^2(s)$} & \{s_0, s_1,\ldots, s_k\} & \{s'_{i} \mid k\leq i <e+1\}  \\ \hline
\textrm{Noyau $N(s)$} & \{s_k\} & \{s'_k\}   \\ \hline 
\end{array}$$

Comme ${\cal N} = {\cal E}(S)$, on a $D(\G) = \{1\}$ et $G = Aut(\G)$ (cf. remarque \ref{remarque finale} b). On a donc $Aut([S],[S]_c) = K^\circ(\G) \rtimes Aut(\G)$.\\

Supposons $e$ fini. On vérifie alors que $Aut(\G) = <\s> \approx \F_2$, où $\s$ envoie $s_k$ sur $s'_{e-k}$ et $s'_k$ sur $s_{e-k}$, pour $0 \leq k \leq e$. On a donc $Aut([S],[S]_c)  = K^\circ(\G) \rtimes <\s> \approx  K^\circ(\G) \rtimes \F_2$. 

On voit que $S$ est d'épaisseur $e$, avec $S_k = \{s'_k,s_{e-k}\}$, pour $0\leq k \leq e$. On a donc, d'après les formules (\ref{3}) et (\ref{4}), $K^\circ(\G) = K^\circ_{S_0}(\G) \rtimes (K^\circ_{S_1}(\G) \rtimes (\cdots \rtimes (K^\circ_{S_{e-2}}(\G) \rtimes K^\circ_{S_{e-1}}(\G))\cdots))$ et, pour $0\leq k \leq e-1$, $K^\circ_{S_k}(\G) \approx [(C^2(s'_k))^\star] \times [(C^2(s_{e-k}))^\star] \approx (\F_2)^{e-k} \times (\F_2)^{e-k} \approx (\F_2)^{2(e-k)}$.\\

Supposons $e = \infty$. Alors $Aut(\G) = \{1\}$, donc $Aut([S],[S]_c) = K^\circ(\G)$. 

On voit que $S$ est d'épaisseur infinie, avec $S_k = \{s'_k\}$, pour tout $k\in \N$, et $S_\infty = \{s_k \mid k\in \N\}$. On ne peut donc pas appliquer le théorème \ref{décomp de KY} à $S$ ; cependant, la proposition \ref{decomp} reste valable : si l'on pose $X = \{s_k \mid k \in \N\}$ et $X' = \{s'_k \mid k \in \N\}$, alors $X$ et $X'$ sont saturées, donc d'après la proposition \ref{decomp}, $K^\circ_X(\G)$ et $K^\circ_{X'}(\G)$ sont distingués dans $K^\circ(\G)$ et $K^\circ(\G) = K^\circ_X(\G) \times K^\circ_{X'}(\G)$. 

Le groupe $K^\circ_X(\G)$ s'identifie à un sous-groupe de $Aut([X]) = GL_{\F_2}([X])$ et, plus précisément, si l'on fixe la base $(\{s_0\},\{s_1\}, \ldots )$ de $[X]$, au sous-groupe de $GL_{\F_2}([X])$ constitué des matrices (infinies à droite) triangulaires supérieures avec des 1 sur la diagonale.

De même, le sous-groupe $K^\circ_{X'}(\G)$ s'identifie à un sous-groupe de $Aut([X']) = GL_{\F_2}([X'])$ et, plus précisément, si l'on fixe la base $(\{s'_0\},\{s'_1\}, \ldots )$ de $[X']$, au sous-groupe de $GL_{\F_2}([X'])$ constitué des matrices (infinies à droite) triangulaires inférieures, avec des 1 sur la diagonale et un nombre fini de $1$ dans chaque colonne, inversibles et telles que l'inverse soit de la même forme.
\end{ex}

\section*{Références.}

\begin{description}
\item[\textmd{[B]}] N. {\textsc{Bourbaki}}. \emph{Groupes et Algèbres de Lie, Chapitres IV-VI.} Hermann, Paris, 1968.
\item[\textmd{[BM]}] P. {\textsc{Bahls}}, M. {\textsc{Mihalik}}. \emph{Reflection Independence In Even Coxeter Groups.} Geom. Ded. à paraître.
\item[\textmd{[BMMN]}] N. {\textsc{Brady}}, J.P. {\textsc{MacCammond}}, B. {\textsc{Mühlherr}}, W.D. {\textsc{Neumann}}. \emph{Rigidity of Coxeter groups and Artin groups.} Geom. Ded. \textbf{94} (2002), 91-109.
\item[\textmd{[H]}] T. {\textsc{Hosaka}}. \emph{Determination up to isomorphism of Right-Angled Coxeter systems.} Proc. Japan Ac. \textbf{79} (2003), 33-35.
\item[\textmd{[M]}] B. {\textsc{Mühlherr}}. \emph{Automorphisms of graph-universal Coxeter Groups.} Journal of Algebra \textbf{200} (1998), 629-649.
\item[\textmd{[R]}] D.G. {\textsc{Radcliffe}}. \emph{Rigidity of Right-Angled Coxeter Groups.} eprint arXiv:math/9901049.
\item[\textmd{[R2]}] D.G. {\textsc{Radcliffe}}. \emph{Unique presentation of Coxeter groups and related groups.} Ph.D. Thesis, University of Wisconsin, Milwaukee (2001).
\item[\textmd{[T]}] J. {\textsc{Tits}}. \emph{Sur le groupe des automorphismes de certains groupes de Coxeter.} Journal of Algebra \textbf{113} (1988), 346-357.
\end{description}

\end{document}